\documentclass[12pt]{amsart}
\usepackage{amscd}
\usepackage{xypic}
\usepackage{color}
\def\ver{s-hammocks-u.tex}

\input prepictex   \input pictex    \input postpictex

\DeclareMathOperator\ind{ind}

\DeclareMathOperator\rad{rad}
\DeclareMathOperator\soc{soc}

\DeclareMathOperator\Ker{Ker}
\DeclareMathOperator\Cok{Cok}
\DeclareMathOperator\Cov{Cov}
\DeclareMathOperator\Hom{Hom}
\renewcommand\hom{{\rm hom}\,}
\DeclareMathOperator\End{End}

\renewcommand\Im{{\rm Im}\,}

\renewcommand\mod{{\rm mod}}

 \def \boldit#1{\textit{\textbf{#1}}}
\swapnumbers

\newtheoremstyle{mytheorems}{9pt}{6pt}{\itshape}{0pt}{\sc}{.}{ }{}
\newtheoremstyle{myremarks}{6pt}{3pt}{\normalfont}{0pt}{\it}{:}{ }{}

\theoremstyle{mytheorems}
\newtheorem{theorem}{Theorem}[section]
\newtheorem{lemma}[theorem]{Lemma}
\newtheorem{claim}[theorem]{Claim}

\newtheorem{proposition}[theorem]{Proposition}

\theoremstyle{myremarks}

\newtheorem{example}[theorem]{Example}
\newtheorem*{remark}{Remark}
\newtheorem*{defin}{Definition}

\newcommand\comment[1]{}
\newenvironment{red}{}{}
\newcommand\mylabel[1]{\label{#1}}

\definecolor{darkgreen}{rgb}{0,0.5,0}
\def\blue#1{{#1}}
\def\green#1{{#1}}

\def\arr#1#2{\arrow <2mm> [0.25,0.75] from #1 to #2 }
\def\sq{\plot 0 0  1 0  1 1  0 1  0 0 / }
\def\bul{$\sssize\bullet$}
\def\sssize{\scriptscriptstyle}
\newcounter{boxsize}
\def\sbullet{\makebox(0,0){$\scriptstyle\bullet$}}
\newcommand\smbox{\put(0,0){\line(1,0){\value{boxsize}}}%
  \put(\value{boxsize},0){\line(0,1){\value{boxsize}}}%
  \put(0,0){\line(0,1){\value{boxsize}}}%
  \put(0,\value{boxsize}){\line(1,0){\value{boxsize}}}}
\def\fac{\multiput{} at 0 -1  0 5 /
  \multiput{\sq} at 0 1  0 2 /
  \put{\bul} at .5 1 }
\def\fae{\multiput{} at 0 -1  0 5 /
  \multiput{\sq} at 0 0  0 1  0 2  0 3 /
  \put{\bul} at .5 1 }
\def\faf{\multiput{} at 0 -1  0 5 /
  \multiput{\sq} at 0 0  0 1  0 2  0 3  1 1  1 2 /
  \multiput{\bul} at .5 2  1.5 2 /
  \plot .5 2  1.5 2 / }
\def\fbb{\multiput{} at 0 -1  0 5 /
  \put{\sq} at   0 2 }
\def\fbd{\multiput{} at 0 -1  0 5 /
  \multiput{\sq} at 0 0  0 1  0 2  0 3  1 1  1 2 /
  \multiput{\bul} at .5 2  1.5 2  1.5 1 /
  \plot .5 2  1.5 2 / }
\def\fbg{\multiput{} at 0 -1  0 5 /
  \multiput{\sq} at  0 1  0 2  0 3  / }
\def\fcc{\multiput{} at 0 -1  0 5 /
  \multiput{\sq} at 0 0  0 1  0 2  0 3   1 2 /
  \multiput{\bul} at .5 2  1.5 2 /
  \plot .5 2  1.5 2 / }
\def\fce{\multiput{} at 0 -1  0 5 /
  \multiput{\sq} at  0 1  0 2   /
  \put{\bul} at .5 2 }
\def\fcf{\multiput{} at 0 -1  0 5 /
  \multiput{\sq} at  0 1  0 2  0 3 /
  \put{\bul} at .5 1 }
\def\fch{\multiput{} at 0 -1  0 5 /
  \multiput{\sq} at   0 1  0 2  0 3  0 4 / }
\def\fdb{\multiput{} at 0 -1  0 5 /
  \multiput{\sq} at 0 0  0 1  0 2  0 3 /
  \put{\bul} at .5 2 }
\def\fdd{\multiput{} at 0 -1  0 5 /
  \multiput{\sq} at  0 1  0 2  0 3   1 2 /
  \multiput{\bul} at .5 2  1.5 2 /
  \plot .5 2  1.5 2 / }
\def\fdg{\multiput{} at 0 -1  0 5 /
  \multiput{\sq} at  0 1  0 2  0 3  0 4 /
  \put{\bul} at .5 1 }
\def\fea{\multiput{} at 0 -1  0 5 /
  \multiput{\sq} at 0 0  0 1  0 2  0 3 /
  \put{\bul} at .5 3 }
\def\fec{\multiput{} at 0 -1  0 5 /
  \multiput{\sq} at   0 1  0 2  0 3 /
  \put{\bul} at .5 2 }
\def\fee{\multiput{} at 0 -1  0 5 /
  \multiput{\sq} at  0 2  0 3 / }
\def\fef{\multiput{} at 0 -1  0 5 /
  \multiput{\sq} at  0 1  0 2  0 3  0 4  1 2 /
  \multiput{\bul} at .5 2  1.5 2 /
  \plot .5 2  1.5 2 / }
\def\ffb{\multiput{} at 0 -1  0 5 /
  \multiput{\sq} at  0 1  0 2  0 3 /
  \put{\bul} at .5 3 }
\def\ffd{\multiput{} at 0 -1  0 5 /
  \multiput{\sq} at   0 1  0 2  0 3  0 4   1 2  1 3  /
  \multiput{\bul} at .5 2  1.5 2 /
  \plot .5 2  1.5 2 / }
\def\ffg{\multiput{} at 0 -1  0 5 /
  \put{\sq} at 0 2
  \put{\bul} at .5 2 }
\def\arqfourblack{%
  \multiput{\scale{\put{\fac}  at .8 0  \put{} at 1 0 }} at  4 4  16 9 /
  \multiput{\scale{\put{\fae}  at .8 0  \put{} at 1 0 }} at  4 7  16 7 /
  \multiput{\scale{\put{\faf}  at 0 0  \put{} at 1 0 }} at  4 9  16 4 /
  \multiput{\scale{\put{\fbb}  at .8 0  \put{} at 1 0 }} at  6 2  18 11 /
  \multiput{\scale{\put{\fbd}  at 0 0  \put{} at 1 0 }} at  6 6  18 6 /
  \multiput{\scale{\put{\fbg}  at .8 0  \put{} at 1 0 }} at  6 11 18 2 /
  \multiput{\scale{\put{\fcc}  at 0 0  \put{} at 1 0 }} at  8 4  20 9 /
  \multiput{\scale{\put{\fce}  at .8 0  \put{} at 1 0 }} at  8 7  20 7 /
  \multiput{\scale{\put{\fcf}  at .8 0  \put{} at 1 0 }} at  8 9  20 4 /
  \multiput{\scale{\put{\fch}  at .8 0  \put{} at 1 0 }} at  8 13  20 0 /
  \multiput{\scale{\put{\fdb}  at .8 0  \put{} at 1 0 }} at  10 2  22 11 /
  \multiput{\scale{\put{\fdd}  at 0 0  \put{} at 1 0 }} at  10 6  22 6 /
  \multiput{\scale{\put{\fdg}  at .8 0  \put{} at 1 0 }} at  10 11 22 2 /
  \multiput{\scale{\put{\fea}  at .8 0  \put{} at 1 0 }} at  0 13  12 0  24 13 /
  \multiput{\scale{\put{\fec}  at .8 0  \put{} at 1 0 }} at  0 9  12 4  24 9 /
  \multiput{\scale{\put{\fee}  at .8 0  \put{} at 1 0 }} at  0 7  12 7  24 7 /
  \multiput{\scale{\put{\fef}  at 0 0  \put{} at 1 0 }} at  0 4  12 9  24 4 /
  \multiput{\scale{\put{\ffb}  at .8 0  \put{} at 1 0 }} at  2 11  14 2 /
  \multiput{\scale{\put{\ffd}  at 0 0  \put{} at 1 0 }} at  2 6  14 6 /
  \multiput{\scale{\put{\ffg}  at .8 0  \put{} at 1 0 }} at  2 2  14 11 /
  \multiput{ \arr {.6 2.4} {1.4 1.6} } at 0 2  0 11  2 4  2 9  4 2  6 4  6 9
    8 2  8 11  10 0  10 4  10 9  12 2  14 4  14 9  16 2  18 0
    18 4  18 9  20 2  22 4  22 9 /
  \multiput{ \arr {.6 1.6} {1.4 2.4} } at 0 4  0 9  2 2  4 4  4 9  6 2  6 11
    8 4  8 9  10 2  12 0  12 4  12 9  14 2  16 4  16 9  18 2  20 0  20 4  20 9
    22 11  22 2 /
  \multiput{ \arr {.6 1.3} {1.4 .8} } at 0 6  4 6  8 6  12 6  16 6  20 6 /
  \multiput{ \arr {.6 .8} {1.4 1.3} } at 2 6  6 6  10 6  14 6  18 6  22 6 /
  \multiput{ \arr {.6 3.6} {1.4 2.6} } at 0 6  4 6  8 6  12 6  16 6  20 6 /
  \multiput{ \arr {.6 2.6} {1.4 3.6} } at 2 6  6 6  10 6  14 6  18 6  22 6 /
  \setdashes<3pt>
  \multiput{ \plot 0 0  0 3 /  \plot 0 5.5  0 6.8 /  \plot 0 10  0 11.8 / } at
  0 5.7  24 5.7 /
  \plot 12.2 1  12.2 3 /  \plot 12.2 5  12.2 6.3 /  \plot 12.2 10.5  12.2 13 /
  \setdots<2pt>
  \multiput{ \plot 0 0  1 0 / } at  0 2  23 2   /
  \multiput{ \plot 0 0  2 0 / } at  3 2  3 11  7 2  11 11  15 2  15 11  19 11
    1 7  5 7  9 7  13 7  17 7  21 7 /
}
\def\arqfourone{%
  \multiput{\scale{\put{\fac}  at .8 0  \put{} at 1 0 }} at  4 4 /
  \multiput{\scale{\put{\fae}  at .8 0  \put{} at 1 0 }} at  4 7 /
  \multiput{\scale{\put{\faf}  at 0 0  \put{} at 1 0 }} at  4 9 /
  \multiput{\scale{\put{\fbb}  at .8 0  \put{} at 1 0 }} at  6 2 /
  \multiput{\scale{\put{\fbd}  at 0 0  \put{} at 1 0 }} at  6 6  /
  \multiput{\scale{\put{\fbg}  at .8 0  \put{} at 1 0 }} at  6 11 /
  \multiput{\scale{\put{\fcc}  at 0 0  \put{} at 1 0 }} at  8 4  /
  \multiput{\scale{\put{\fce}  at .8 0  \put{} at 1 0 }} at  8 7 /
  \multiput{\scale{\put{\fcf}  at .8 0  \put{} at 1 0 }} at  8 9 /
  \multiput{\scale{\put{\fch}  at .8 0  \put{} at 1 0 }} at  8 13 /
  \multiput{\scale{\put{\fdb}  at .8 0  \put{} at 1 0 }} at  10 2  -2 11 /
  \multiput{\scale{\put{\fdd}  at 0 0  \put{} at 1 0 }} at  10 6  -2 6 /
  \multiput{\scale{\put{\fdg}  at .8 0  \put{} at 1 0 }} at  10 11 -2 2 /
  \multiput{\scale{\put{\fea}  at .8 0  \put{} at 1 0 }} at  0 13 /
  \multiput{\scale{\put{\fec}  at .8 0  \put{} at 1 0 }} at  0 9  /
  \multiput{\scale{\put{\fee}  at .8 0  \put{} at 1 0 }} at  0 7  /
  \multiput{\scale{\put{\fef}  at 0 0  \put{} at 1 0 }} at  0 4  /
  \multiput{\scale{\put{\ffb}  at .8 0  \put{} at 1 0 }} at  2 11 /
  \multiput{\scale{\put{\ffd}  at 0 0  \put{} at 1 0 }} at  2 6  /
  \multiput{\scale{\put{\ffg}  at .8 0  \put{} at 1 0 }} at  2 2 /
  \multiput{ \arr {.6 2.4} {1.4 1.6} } at 0 2  0 11  2 4  2 9  4 2  6 4  6 9
    8 2  8 11  -2 4  -2 9 /
  \multiput{ \arr {.6 1.6} {1.4 2.4} } at 0 4  0 9  2 2  4 4  4 9  6 2  6 11
    8 4  8 9  -2 2  -2 11 / 
  \multiput{ \arr {.6 1.3} {1.4 .8} } at 0 6  4 6  8 6 /
  \multiput{ \arr {.6 .8} {1.4 1.3} } at 2 6  6 6  -2 6 /
  \multiput{ \arr {.6 3.6} {1.4 2.6} } at 0 6  4 6  8 6 /
  \multiput{ \arr {.6 2.6} {1.4 3.6} } at 2 6  6 6  -2 6 /
  \setdashes<3pt>
  \plot -1.85 0.3  -1.85 1.2 /  \plot -1.85 3.6  -1.85 4.8 /  \plot -1.85 7  -1.85 9.5 /  \plot -1.85 12  -1.85 14.2 /
  \plot 10.15 0.3  10.15 .6 /  \plot 10.15 3.2  10.15 5 /  \plot 10.15 7  10.15 10 /  \plot 10.15 12.5  10.15 14.2 /
  \setdots<2pt>
  \multiput{ \plot 0 0  1 0 / } at  9 7  -2 7 /
  \multiput{ \plot 0 0  2 0 / } at  -1 2  3 2  3 11  7 2
    1 7  5 7   /
}

\def\emod#1#2#3#4#5#6{$\scriptscriptstyle {#1\,#3\,#6\atop\,#2#4#5\,}$}
\def\arqesix{ %
\multiput{} at -1 -1  12 9 /
\put{\emod000010} at 0 2
\put{\emod000100} at 0 3
\put{\emod010000} at 0 6
\put{\emod000011} at 1 0
\put{\emod011110} at 1 4
\put{\emod110000} at 1 8
\put{\emod011111} at 2 2
\put{\emod011010} at 2 3
\put{\emod111110} at 2 6
\put{\emod011100} at 3 0
\put{\emod122121} at 3 4
\put{\emod001110} at 3 8
\put{\emod122110} at 4 2
\put{\emod111111} at 4 3
\put{\emod012121} at 4 6
\put{\emod111010} at 5 0
\put{\emod123221} at 5 4
\put{\emod011011} at 5 8
\put{\emod112121} at 6 2
\put{\emod012110} at 6 3
\put{\emod122111} at 6 6
\put{\emod001111} at 7 0
\put{\emod123121} at 7 4
\put{\emod111100} at 7 8
\put{\emod012111} at 8 2
\put{\emod111011} at 8 3
\put{\emod112110} at 8 6
\put{\emod011000} at 9 0
\put{\emod112111} at 9 4
\put{\emod001010} at 9 8
\put{\emod111000} at 10 2
\put{\emod001100} at 10 3
\put{\emod001011} at 10 6
\put{\emod100000} at 11 0
\put{\emod001000} at 11 4
\put{\emod000001} at 11 8
\arr{.3 1.4}{.7 .6}
\arr{2.3 1.4}{2.7 .6}
\arr{4.3 1.4}{4.7 .6}
\arr{6.3 1.4}{6.7 .6}
\arr{8.3 1.4}{8.7 .6}
\arr{10.3 1.4}{10.7 .6}
\arr{.3 2.6}{.7 3.4}
\arr{2.3 2.6}{2.7 3.4}
\arr{4.3 2.6}{4.7 3.4}
\arr{6.3 2.6}{6.7 3.4}
\arr{8.3 2.6}{8.7 3.4}
\arr{10.3 2.6}{10.7 3.4}
\arr{.3 3.3}{.7 3.7}
\arr{2.3 3.3}{2.7 3.7}
\arr{4.3 3.3}{4.7 3.7}
\arr{6.3 3.3}{6.7 3.7}
\arr{8.3 3.3}{8.7 3.7}
\arr{10.3 3.3}{10.7 3.7}
\arr{.3 5.4}{.7 4.6}
\arr{2.3 5.4}{2.7 4.6}
\arr{4.3 5.4}{4.7 4.6}
\arr{6.3 5.4}{6.7 4.6}
\arr{8.3 5.4}{8.7 4.6}
\arr{10.3 5.4}{10.7 4.6}
\arr{.3 6.6}{.7 7.4}
\arr{2.3 6.6}{2.7 7.4}
\arr{4.3 6.6}{4.7 7.4}
\arr{6.3 6.6}{6.7 7.4}
\arr{8.3 6.6}{8.7 7.4}
\arr{10.3 6.6}{10.7 7.4}
\arr{1.3 .6}{1.7 1.4}
\arr{3.3 .6}{3.7 1.4}
\arr{5.3 .6}{5.7 1.4}
\arr{7.3 .6}{7.7 1.4}
\arr{9.3 .6}{9.7 1.4}
\arr{1.3 3.4}{1.7 2.6}
\arr{3.3 3.4}{3.7 2.6}
\arr{5.3 3.4}{5.7 2.6}
\arr{7.3 3.4}{7.7 2.6}
\arr{9.3 3.4}{9.7 2.6}
\arr{1.3 3.7}{1.7 3.3}
\arr{3.3 3.7}{3.7 3.3}
\arr{5.3 3.7}{5.7 3.3}
\arr{7.3 3.7}{7.7 3.3}
\arr{9.3 3.7}{9.7 3.3}
\arr{1.3 4.6}{1.7 5.4}
\arr{3.3 4.6}{3.7 5.4}
\arr{5.3 4.6}{5.7 5.4}
\arr{7.3 4.6}{7.7 5.4}
\arr{9.3 4.6}{9.7 5.4}
\arr{1.3 7.4}{1.7 6.6}
\arr{3.3 7.4}{3.7 6.6}
\arr{5.3 7.4}{5.7 6.6}
\arr{7.3 7.4}{7.7 6.6}
\arr{9.3 7.4}{9.7 6.6}
\setdots<2pt>
\plot 1.3 0  2.7 0 /
\plot 3.3 0  4.7 0 /
\plot 5.3 0  6.7 0 /
\plot 7.3 0  8.7 0 /
\plot 9.3 0  10.7 0 /
\plot 1.3 8  2.7 8 /
\plot 3.3 8  4.7 8 /
\plot 5.3 8  6.7 8 /
\plot 7.3 8  8.7 8 /
\plot 9.3 8  10.7 8 /
\plot 0.3 3  1.7 3 /
\plot 2.3 3  3.7 3 /
\plot 4.3 3  5.7 3 /
\plot 6.3 3  7.7 3 /
\plot 8.3 3  9.7 3 /
}

\begin{document}
{\footnotesize [\ver, \today]}

        \vglue1truecm
        \begin{center}{\Large Hammocks to Visualize\\[1ex]
            the Support of Finitely Presented Functors}

        \bigskip
        \normalsize Markus Schmidmeier

\bigskip
\parbox{11cm}{%
  \small {\it Abstract:}  Many properties of a module can be expressed in terms of the dimension of
  the vector space obtained by applying a finitely presented functor to that module.
  For example, the dimension of the kernel, image 
  or cokernel of the multiplication map given by an algebra element;
  or the number of summands of a certain type
  when the module is considered a module over a subalgebra.
  When the indecomposable modules over the algebra are arranged in the Auslander-Reiten quiver,
  the support of the finitely presented functor typically has the shape of a hammock, spanned
  between sources and sinks.  There may also be tangents which are meshes where the
  hammock function at the middle term exceeds the sum of the values at the start and end terms.
  We describe how sources, sinks and tangents of the hammock relate to the modules which define
  the projective resolution of the finitely presented functor.
  The key tool is the Cokernel Complex Lemma which links the values of the hammock function
  to the Auslander-Reiten structure of the category.
  We are also interested in exact subcategories of module categories which have Auslander-Reiten
  sequences.  Our examples include quiver representations and invariant subspaces of
  nilpotent linear operators.

  \medskip {\it Keywords:} Hammock, Auslander-Reiten quiver, finitely presented functor,
  invariant subspace.

  \medskip {\it MSC 2020:} 16G70 (primary), 16E05, 47A15
}        
        \end{center}
        \bigskip

\section{Introduction}

Let $\Lambda$ be a finite dimensional algebra over a field $k$.
Our interest is in the category $\mod\Lambda$
of finite dimensional (right) $\Lambda$-modules, or more general, in 
a full exact extension-closed subcategory $\mathcal C$ of $\mod\Lambda$ which,
in particular, has Auslander-Reiten sequences.
Thus, we can picture the category $\mathcal C$ as a directed graph,
the Auslander-Reiten quiver $\Gamma=(\Gamma_0,\Gamma_1,\tau)$, which comprises the
information how homomorphisms in $\mathcal C$ factor locally.
The vertex set $\Gamma_0=\ind\mathcal C$ consists of the isomorphism classes
of indecomposable objects, the arrows in $\Gamma_1$ represent irreducible morphisms,
and the translation $\tau$ is a partial map on $\Gamma_0$. 
Each Auslander-Reiten sequence
$\mathcal A: 0\to A \to \bigoplus B_i\to C\to 0$
contributes a mesh in $\Gamma$ in the sense that $A=\tau C$ and
the irreducible morphisms starting at $A$ correspond to irreducible morphisms ending in $C$.
For more details and examples see \cite[Chapter~VII]{ARS}.

\medskip
Our aim is to visualize the support of a finitely presented functor $E:\mod\Lambda\to \mod k$
as a subset of $\Gamma_0$.  Such a functor is given by 
an exact sequence of Hom-functors,
$$0\longrightarrow (Z,-)\longrightarrow (Y,-)\longrightarrow (X,-)\longrightarrow E
\longrightarrow 0,$$
where $X,Y,Z\in\mod\Lambda$ are finitely generated modules
and where we write $(X,-)$ as abbreviation for $\Hom_\Lambda(X,-)$.
The sequence is induced by a short right exact sequence
$$\mathcal E: \quad X\stackrel u\longrightarrow Y\longrightarrow Z\longrightarrow 0$$
of $\Lambda$-modules.  

\begin{example}
  \mylabel{ex-intro}
  Before we state our main result, we illustrate this aim in an example.
  Let $\Lambda$ be the path algebra of the quiver $\mathbb E_6$ in bipartite orientation as in Section~\ref{sec-cok}.
  Then the category of $\Lambda$-modules is equivalent to the category of systems $V=(V_i,V_\alpha)$
  consisting of six vector spaces and five linear maps, arranged as follows.
  $$
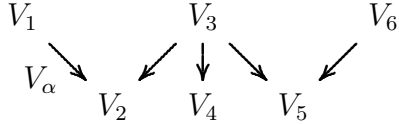
\beginpicture\setcoordinatesystem units <1.2cm,1.2cm>
  \put{$V_1$} at 0 1
  \put{$V_2$} at 1 0
  \put{$V_3$} at 2 1
  \put{$V_4$} at 2 0
  \put{$V_5$} at 3 0
  \put{$V_6$} at 4 1
  \arr{0.3 0.7}{0.7 0.3}
  \arr{1.7 0.7}{1.3 0.3}
  \arr{2 0.7}{2 0.3}
  \arr{2.3 0.7}{2.7 0.3}
  \arr{3.7 0.7}{3.3 0.3}
  \put{$V_\alpha$} at 0.2 0.3 
  \endpicture$$
  It turns out that
  there are 36 indecomposable $\Lambda$-modules, each is determined uniquely by its dimension
  type $(\dim V_i)_i$; for example, the indecomposable module $X$ which is zero in each component
  except the second where $X_2=k$ has dimension type \emod010000.
  In the Auslander-Reiten quiver in Figure~\ref{fig-intro}, the indecomposable objects
  are represented by their dimension types, the irreducible morphisms by arrows,
  and the translation $\tau$ maps each non-projective object $C$ to the object $A$
  on the left so that the Auslander-Reiten sequence $0\to A\to \bigoplus B_i\to C\to 0$
  gives rise to a mesh in the graph.  

  \smallskip
  In this example
  we are interested in the (finitely presented) functor $E:\mod\Lambda\to\mod k$ which assigns to
  each object $V$ the cokernel $\Cok V_\alpha$ of the leftmost linear map.  The collection of modules
  $V$ which satisfy $E(V)\neq 0$ is encircled.

  \begin{figure}[ht]
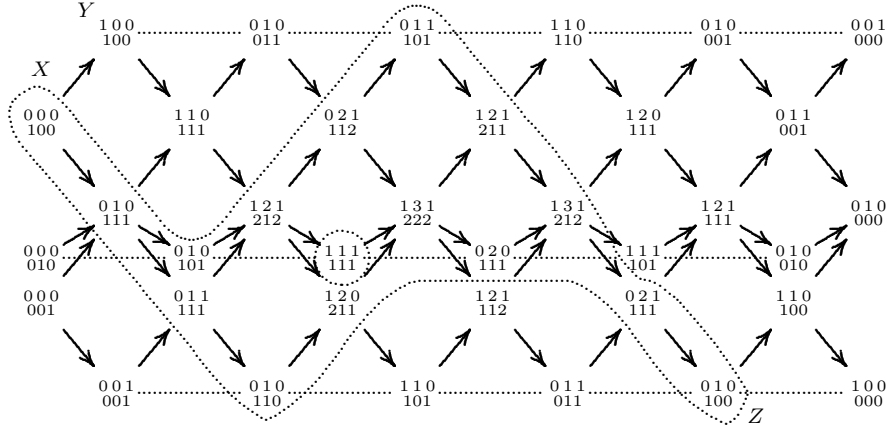

$$
\beginpicture\setcoordinatesystem units <1cm,.6cm>
\arqesix
\put{$\scriptstyle X$} at 0 7.2
\put{$\scriptstyle Y$} at .6 8.5
\put{$\scriptstyle Z$} at 9.5 -.5
\setquadratic
\setdots<2pt>
\circulararc 360 degrees from 4 2.4 center at 4 3  
\plot 1 3  1.5 2  2 1  2.5 0  3 -.6  3.5 0  4 1  4.5 2  5 2.5  5.5 2.5  6 2.5  6.5 2.5
7 2.5 7.5 2  8 1  8.5 0  9 -.6  9.2 -.6  9.4 0  9 1  8.5 2  8.2 2.4  8 2.5  7.7 3
7.5 3.7  7 5  6.5 6  6 7  5.5 8  5 8.6  4.5 8  4 7 3.5 6  3 5  2.5 4  2 3.4  1.5 4  1 5
0.5 6  0.1 6.7  -.2 6.7  -.4 6.4  -.4 6  -.3 5.6  0 5  0.5 4  1 3 /
\endpicture
$$
    \caption{The hammock for the functor $E:V\mapsto \Cok V_\alpha$}
    \mylabel{fig-intro}
\end{figure}

  \smallskip
The encircled region is a hammock:  It has a source at $X$,
similarly, there is a sink at the module $Z'$ of dimension type \emod011000, and we will see
that the functor $E$ is additive on each mesh except the ones involving $X$, $Y$ and $Z'$.
We call the mesh ending at the module marked $Y$ a tangent for $E$.
Denote by $Z=\tau^{-1}Z'$ the module on the right hand side of $Z'$, of dimension type \emod 100000.

\smallskip
Our main result relates sinks, sources and tangents of the hammock with the modules defining
the presentation of the functor.  In fact, in this case the functor $E$ has projective presentation
$$0\longrightarrow (Z,-)\longrightarrow (Y,-)\longrightarrow (X,-)\longrightarrow E\longrightarrow 0.$$
\end{example}

\begin{defin}Given a finitely presented functor $E:\mod\Lambda\to \mod k$, we call the map
$$e_*: \quad \Gamma_0\to \mathbb N_0, \quad M\mapsto \dim EM,$$
the \boldit{hammock function} on $\Gamma$ (or on $\mathcal C$)
for the functor $E$ (or for the sequence $\mathcal E$).  A \boldit{hammock}
is the support of a hammock function.
  \begin{enumerate}
  \item A vertex $C\in\Gamma_0$ is an \boldit{(isolated) source} for $e_*$ if $e_*(C)>\sum_{B\in C^-} e_*(B)$.
  \item A vertex $A\in\Gamma_0$ is an \boldit{(isolated) sink} for $e_*$ if $e_*(A)>\sum_{B\in A^+} e_*(B)$.
  \item A mesh $\mathcal M:A\to \bigoplus B_i\to C$ in $\Gamma$ is an \boldit{(isolated) tangent} for $e_*$
    if $\sum e_*(B_i)>e_*(A)+e_*(C)$.
  \end{enumerate}
\end{defin}

\medskip
For $A,C\in \Gamma_0$,
we denote by $A^+$ and $C^-$ the multisets of successors
and predecessors of $A$ and $C$ with respect to irreducible morphisms in $\Gamma_1$,
respectively.
A \boldit{mesh} is a sequence in $\Gamma$ that has the form
$0\to \bigoplus_{B\in C^-}B\to C$ if $C$ is an indecomposable projective object in $\mathcal C$,
or $A\to \bigoplus_{B\in A^+}B\to 0$ if $A$ is indecomposable injective.
Otherwise, the mesh $\tau C\to \bigoplus_{B\in C^-}B\to C$ is given by an Auslander-Reiten sequence.

\medskip
Here is our main result about hammock functions.

\begin{theorem}
  \mylabel{theorem-isolated}
  Let $\mathcal E: X\stackrel u\to Y\to Z\to 0$ be a short right exact sequence
  and $e_*$ the hammock function for $\mathcal E$ on $\Gamma$.
  \begin{enumerate}
  \item If $C$ is an isolated source for $e_*$ then $C$ is a direct summand of $X$.
  \item If $A$ is an isolated sink for $e_*$ then $A$ is an injective object in
    $\mathcal C$ or $C=\tau^{-1}A$ is a direct summand of $Z$.
  \item If $\mathcal M:A\to \bigoplus B_i\to C$ is an isolated tangent for $e_*$
    then either $A$ is indecomposable injective or $C$ is
    an indecomposable direct summand of $Y$.
  \item If $\mathcal M:A\to \bigoplus B_i\to C$ is a mesh in $\Gamma$
    with $C\neq 0$ such that
    $C$ does not occur as a direct summand of $X\oplus Y\oplus Z$ then
    the hammock function $e_*$ is additive on $\mathcal M$
    in the sense that  $e_*(A)+e_*(C)=\sum e_*(B_i)$ holds. 
  \end{enumerate}
\end{theorem}

\medskip
The theorem describes how the terms of $\mathcal E:X\to Y\to Z\to 0$ control the shape of the hammock
given by the support of the functor $E$ on $\Gamma$.
The hammock is clearly visible
if $\mathcal C$ is representation-directed (so there are no oriented cycles in $\Gamma$):

\medskip
Tracing the hammock function $e_*$ on $\Gamma$ from left to right, the source or sources of the hammock
occur at the indecomposable summands of $X$ (1).  The values of $e_*$ can be computed successively by using the
mesh relation (4), except when summands of $X\oplus Y\oplus Z$ are encountered.  Summands of $Y$ are
endpoints of meshes where a middle term is touched by the hammock (3).  Finally, the hammock terminates
at one or several vertices which are injective or of the form $\tau Z'$ where $Z'$ is a summand of $Z$ (2).

\medskip
We use the prefix {\it isolated} since indecomposable objects may occur as direct
summands in several of the terms $X$, $Y$, $Z$ in the resolution for $E$.
Such sources, sinks and tangents may not be
identifiable by the numerical condition in the definition alone.  See Section~\ref{sec-no-isolated} for
an example.

\medskip
Theorem~\ref{theorem-isolated} deals with hammocks which visualize the support
of a finitely presented covariant functor.  The corresponding dual result
for finitely presented contravariant functors is stated as Theorem~\ref{theorem-isolated-dual}.

\bigskip
There is a different way to describe the modules in the hammock.

\smallskip
Suppose $\mathcal E:0\to X\stackrel u\to Y\to Z\to 0$ is given by a non-split short exact sequence
with $Z$ indecomposable non-projective.  Then
the  modules $M\in\mathcal C$ which occur in  the hammock given by $\mathcal E$
(which hence satisfy by definition that $0\neq\Cok\Hom(u,M):\Hom(Y,M)\to\Hom(X,M)$) are characterized by lying
between one of the summands of $X$ and
of $\tau Z$ in the sense that there is a map $f:X\to M$ which does not factor over $u$
and a map $g:M\to \tau Z$ such that the product $gf:X\to \tau Z$ does not
factor over $u$.  Our result generalizes, but is weaker than, the corresponding statement
\cite[Corollary~5]{rv}.

\begin{proposition}
  \mylabel{prop-hammock-factoring}
  Let $\mathcal E: 0\to X\stackrel u\to Y\to Z\to 0$
  be a non-split short exact sequence in $\mathcal C$ with $Z$
  indecomposable. Then for any $M\in\mathcal C$, the bilinear form given by composition
  $$\Hom(M,\tau Z)\times\Cok(u,M) \longrightarrow \Cok(u,\tau Z)$$
  is right non-degenerate.
\end{proposition}

The dual result for finitely presented contravariant functors is stated as
Proposition~\ref{prop-hammock-factoring-dual}.

\medskip
In the proof of Theorem~\ref{theorem-isolated}, we relate the Auslander-Reiten structure
of the category $\mathcal C$ with the modules in the projective resolution of the functor $E$.
The following result is our key tool.  It will be shown in Section~\ref{sec-snake}.

\begin{proposition}[Cokernel Complex Lemma]\mylabel{proposition-isomorphic-homology}
Consider the commutative diagram with exact rows and columns.
$$
\xymatrix{ & & & 0 \ar[d] \\
  & & A \ar[r]^\alpha \ar[d] & A''\ar[d] \\
  & B' \ar[r] \ar[d]_{g'} & B \ar[r]^\beta \ar[d]_g & B'' \ar[d]_{g''} \\
  0 \ar[r] & C' \ar[r] & C \ar[r]^\gamma & C'' }
$$

The complexes $0\to \Cok\alpha \to \Cok\beta \to \Cok\gamma \to 0$ and
$0\to \Cok g' \to \Cok g \to \Cok g'' \to 0$ given by taking
the cokernels of the horizontal maps
at the right and the vertical maps at the bottom, respectively, have in each
position isomorphic homology.
\end{proposition}

\bigskip
The second part of this paper consists of examples.  We consider a variety of quantities
related to modules which can be measured by hammock functions.

\medskip
In Section~\ref{sec-quiver} we consider multiplication maps.
For an element $a$ in an algebra $\Lambda$ and $M$ a finite dimensional $\Lambda$-module,
the multiplication by $a$ gives rise to a linear map $\mu_a:M\to M$.  
We describe the hammock functions given by the kernel, $\dim\Ker\mu_a$,
the image, $\dim\Im\mu_a$, and the cokernel, $\dim\Cok\mu_a$.
If $\Lambda=kQ$ is the path algebra of a quiver, $i\in Q_0$ a vertex, $\alpha\in Q_1$ an arrow
and $M=((M_i)_{i\in Q_0},(M_\alpha)_{\alpha\in Q_1})$ a representation,
then quantities like $\dim M_i$, $\dim \Ker M_\alpha$,
$\dim\Im M_\alpha$, $\dim\Cok M_\alpha$ give rise
to hammock functions; some of them are pictured.
Suppose $p,q$ are composable paths in the quiver, say meeting at the vertex $i$,
such that $M_p\circ M_q=0$.
Then also the homology $\dim\Ker M_p/\Im M_q$ (which is a subspace of $M_i$)
defines a hammock function.

\medskip
In Section~\ref{sec-invariant} we deal with invariant subspaces of nilpotent
linear operators; they are pairs $(V,U)$ where the linear operator $V$ is a
module over some bounded polynomial ring $k[T]/(T^n)$ and the
invariant subspace $U$ is a $k[T]$-submodule of $V$.  Such systems form a full exact subcategory
of a module category which has Auslander-Reiten sequences.  Note that if $n$ is large enough,
then the $k[T]$-modules $U$, $V$, $V/U$ may have many indecomposable direct summands,
even if the pair $(V,U)$ itself is indecomposable.  For each $1\leq m< n$,
we exhibit a hammock function which counts the number of summands in $V$ (and similarly
in $U$, $V/U$) which are isomorphic to $k[T]/(T^m)$.
For $n=4$ we picture some such functions as hammocks in the Auslander-Reiten quiver.
Our last example is a hammock for which the source is not isolated.

\medskip
\begin{remark}
  Hammocks have been introduced by Sheila Brenner \cite{b} to study,
  for a representation directed algebra, the collection of indecomposable modules
  which contain a given simple module as a composition factor.
  Regarding finitely presented functors, there has been interest in their support, 
  see for example the paper \cite{ar-uniserial} on uniserial functors.
  The link to poset representations is established in \cite{rv}; a decomposition
  of the poset given by the hammock is studied in \cite{sc}.
  The minimal elements in hammocks are described in \cite{xi}.
  Y.\ Lin \cite{lin1,lin2} uses hammocks to trace algorithms
  introduced by Nazarova and Roiter, and by Zavadskii.  In \cite{brue},
  kit algebras are introduced as representation-directed algebras for which
  each hammock given by the modules which contain a certain simple composition
  factor is a garland.
\end{remark}

\bigskip
We briefly describe the contents of this paper.

\medskip
The Cokernel Complex Lemma is shown in 
Section~\ref{sec-snake}.
For the proof we use relations as they may be occur in a proof of the Snake Lemma.

\medskip
In Section~\ref{sec-hammock} we show the above results about the modules in the hammock
and state in Theorem~\ref{theorem-isolated-dual} and Proposition~\ref{prop-hammock-factoring-dual}
the dual versions which deal with finitely presented contravariant functors.
Our first example are hammocks given by Auslander-Reiten sequences, they are characterized by
having only one vertex.

\medskip
In the last two sections we present applications which show that hammock functions
exhibit a variety of meaningful properties of modules.

\section{The Cokernel Complex Lemma}
\mylabel{sec-snake}

Our main tool to study hammocks is the Cokernel Complex Lemma, stated in the
introduction as Proposition~\ref{proposition-isomorphic-homology}.
For its proof we introduce relations following \cite[Section~26]{rs}.

\smallskip
Recall that for modules $X,Y$, a \boldit{relation} on $Y\times X$ is a submodule
of $Y\oplus X$.
In particular, a map $f:X\to Y$ gives rise to relations $f=\{(f(x),x)|x\in X\}\subset Y\oplus X$
and $f^{-1}=\{(x,f(x))|x\in X\}\subset X\oplus Y$.
Given two relations $u\subset Y\oplus X$ and $v\subset Z\oplus Y$,
the composition $v\circ u$ is the relation given by the submodule
$$v \circ u\;=\; \{ (z,x) \;|\; (y,x)\in u, (z,y)\in v \;\text{for some}\; y\in Y\}
                \;\subset \; Z\oplus X.$$
Clearly, a relation $u$ on $Y\times X$ is given by a map if for every $x\in X$ there
is a unique $y\in Y$ with $(y,x)\in u$; this map is an isomorphism if for every $y\in Y$
there is a unique $x\in X$ with $(y,x)\in u$.

\bigskip
The proof of the following lemma is straightforward and easy.

\begin{lemma}
Consider the diagram with one exact row and one exact column.
$$
\xymatrix{ & A \ar[d]^f & \\
  B' \ar[r]^{\beta'} & B \ar[r]^\beta \ar[d]^g & B'' \\
  & C & }
$$
\begin{enumerate}
\item Let $u\subset C\oplus B''$ be the relation given by $g\circ \beta^{-1}$.
  For $(c,b'')\in u$ we have $b''\in \Im \beta f$ if and only if 
  $c\in\Im g\beta'$.
\item Let $v\subset B'\oplus A$ be the relation given by ${\beta'}^{-1}\circ f$.
  \begin{enumerate}
    \item   Let $a\in A$.  Then $a\in\Ker \beta f$ if and only if there is 
      $b'\in B'$ with $(b',a)\in v$.
    \item  Given $b'\in B'$, we have $b'\in\Ker g \beta'$ if and only if there is 
      $a\in A$ with $(b',a)\in v$. \qed
  \end{enumerate}
\end{enumerate}
\end{lemma}

\medskip
The Cokernel Complex Lemma is a consequence of the following three claims.
In each, we consider parts of the following commutative diagram in which all
rows and columns are exact except the rightmost column and the row at the bottom
(which may be proper complexes).
$$
\xymatrix{ & & & 0 \ar[d] & & \\
  & & A\ar[r]^\alpha \ar[d]_f & A'' \ar[r]^{\alpha''} \ar[d]_{f''} & \tilde A \ar[r] \ar[d]_{\tilde f} & 0 \\
  & B'\ar[r]^{\beta'} \ar[d]_{g'} & B \ar[r]^\beta \ar[d]_g & B'' \ar[r]^{\beta''} \ar[d]_{g''} & \tilde B \ar[r] \ar[d]_{\tilde g} & 0 \\
  0 \ar[r] & C' \ar[r]^{\gamma'} \ar[d]_{h'} & C \ar[r]^\gamma \ar[d]_h & C''\ar[r]^{\gamma''} \ar[d]_{h''} & \tilde C \ar[r] & 0 \\
  & D' \ar[r]^{\delta'} \ar[d] & D \ar[r]^\delta \ar[d] & D'' \ar[d] & & \\
  & 0 & 0 & 0 & & }
$$

\begin{claim} The relation $w=h''\circ{\gamma''}^{-1}\subset D''\times \tilde C$ as in the
  diagram yields an isomorphism $\tilde C/\Im\tilde g\cong D''/\Im\delta$.
  $$
  \xymatrix{ & B'' \ar[r]^{\beta''} \ar[d]_{g''} & \tilde B \ar[r] \ar[d]^{\tilde g} & 0 \\
    C \ar[r]^\gamma \ar[d]_h & C'' \ar[r]^{\gamma''} \ar[d]^{h''} & \tilde C \ar[r] & 0 \\
    D \ar[r]^\delta \ar[d] & D'' \ar[d] & & \\
    0 & 0 & & }
  $$
\end{claim}

\begin{proof}
The relation $w$ has the following properties:

\begin{enumerate}
\item For every $d''\in D''$ there is 
          $\tilde c\in \tilde C$ such that $(d'',\tilde c)\in w$.
\item For every $\tilde c\in \tilde C$ there is $d''\in D''$ 
  such that $(d'',\tilde c)\in w$.
\item For $(d'',\tilde c)\in w$, $d''\in\Im\delta$ if and only if $\tilde c\in\Im \tilde g$.
  Namely, $\tilde c\in\Im \tilde g$ is equivalent to $\tilde c\in\Im \tilde g\beta''=\Im \gamma''g''$
  since $\beta''$ is onto.  By the Lemma, Part (1), this is equivalent to
  $d''\in\Im h''\gamma=\Im\delta h$ and hence to $d''\in\Im\delta$.
\end{enumerate}

We obtain from (2) and (3) that there is a map $\tilde C\to D''/\Im\delta$
which clearly is a homomorphism.  By (1), this 
map is onto and by (3), the map has  kernel $\Im\tilde g$, and hence 
induces an isomorphism $\tilde C/\Im\tilde g\to D''/\Im\delta$.
\end{proof}

\begin{claim}
  The relation $v=h\circ\gamma^{-1}\circ g''\circ {\beta''}^{-1}\subset D\times \tilde B$
  as in the diagram gives rise to an isomorphism
  $\Ker\tilde g/\Im\tilde f\cong \Ker \delta/\Im \delta'$.
 $$
 \xymatrix{ & & A'' \ar[r]^{\alpha''} \ar[d]_{f''} & \tilde A \ar[d]^{\tilde f} \ar[r] & 0 \\
   & B \ar[r]^\beta \ar[d]_g & B'' \ar[r]^{\beta''} \ar[d]_{g''} & \tilde B \ar[d]^{\tilde g} \ar[r] & 0 \\
   C' \ar[r]^{\gamma'} \ar[d]_{h'} & C \ar[r]^\gamma \ar[d]_h & C''\ar[r]^{\gamma''} \ar[d]_{h''} & \tilde C & \\
   D' \ar[r]^{\delta'} \ar[d] & D\ar[r]^\delta \ar[d] & D'' & & \\
   0 & 0 & & & }
 $$
\end{claim}

\begin{proof}
The relation $v$ has the following properties.

\begin{enumerate}
\item Suppose $\tilde b\in \Ker \tilde g$.  There is $d\in \Ker\delta$
          with $(d,\tilde b)\in v$.  Namely, since $\beta''$ is onto, there
          is $b''$ with $\beta''(b'')=\tilde b$. We apply the Lemma, Part (2a),
          center at $C''$,
          to the relation $\gamma^{-1}\circ g''$:
          Since $0=\tilde g\beta''(b'')=\gamma''g''(b'')$, there is $c\in C$ with
          $(c,b'')$ in the relation given by $\gamma^{-1}\circ g''$ 
          and $h''\gamma(c)=0$. Then $(h(c),\tilde b)\in v$ and $h(c)\in \Ker\delta$.
\item Similarly, 
  for $d\in\Ker\delta$ there is $\tilde b\in\Ker \tilde g$ 
  with $(d,\tilde b)\in v$.  
\item For $(d,\tilde b)\in v$, $d\in\Im\delta'$ if and only if $\tilde b\in\Im \tilde f$:
  Let $c''\in C''$ such that $(d,c'')\in h\circ\gamma^{-1}$ and $(c'',\tilde b)\in g''\circ(\beta'')^{-1}$.
  Then $\tilde b\in\Im\tilde f$ is equivalent to 
  $\tilde b\in\Im\tilde f\alpha''=\Im\beta''f''$, since $\alpha''$ is onto.
  By (1) in the Lemma, center at $B''$, this is
  equivalent to $c''\in\Im g''\beta=\Im\gamma g$, which by (1) in the
  Lemma, center at $C$, is
  equivalent to $d\in\Im h\gamma'=\Im\delta'h'$ and hence to $d\in\Im\delta'$.
\end{enumerate}

\smallskip
From (1) and (3) we obtain that there is a map $\Ker\tilde g\to \Ker\delta/\Im\delta'$;
this map is onto by (2) and has kernel $\Im\tilde f$ by (3).
Hence it induces the desired isomorphism
$\Ker\tilde g/\Im\tilde f\cong \Ker\delta/\Im\delta'$.
\end{proof}

\begin{claim}
The relation $u=h'\circ {\gamma'}^{-1}\circ g\circ\beta^{-1}\circ f''\circ{\alpha''}^{-1}$
as in the diagram induces an isomorphism $\Ker\tilde f\to \Ker\delta'$.
$$
\xymatrix{& & & 0 \ar[d] & & \\
  & & A \ar[r]^\alpha \ar[d]_f & A'' \ar[r]^{\alpha''} \ar[d]^{f''} & \tilde A \ar[d]^{\tilde f} \ar[r] & 0 \\
  & B' \ar[r]^{\beta'} \ar[d]_{g'} & B \ar[r]^\beta \ar[d]_g & B'' \ar[r]^{\beta''} \ar[d]^{g''} & \tilde B & \\
 0 \ar[r] & C' \ar[r]^{\gamma'} \ar[d]_{h'} & C \ar[r]^\gamma \ar[d]_h & C'' & & \\
 & D' \ar[r]^{\delta'} \ar[d] & D & & & \\
 & 0 & & & & }
$$
\end{claim}

\begin{proof}
Using the Lemma we obtain the following properties:
\begin{enumerate}
\item For $\tilde a\in\Ker \tilde f$ there is $d'\in \Ker\delta'$ with
  $(d',\tilde a)\in u$:
  Let $\tilde a \in \Ker \tilde f$.
  Since $\alpha''$ is onto, there is $a''\in A''$ with $\alpha''(a'')=\tilde a$.
  By Part (2a), since $\beta''f''(a'')=0$, there is $b\in B$ with $b\in \Ker g''\beta$
  and $\beta(b)=f''(a'')$.
  Again by Part (2a), since $b\in\Ker\gamma g=\Ker g''\beta$, there is  $c'\in C'$ with
  $c'\in \Ker h\gamma'$ and $\gamma'(c')=g(b)$.  Then $h'(c)\in\Ker \delta'$ and
  $(h'(c),\tilde a)\in u$.
\item For $d'\in\Ker \delta'$ there is $\tilde a\in\Ker \tilde f$ with
          $(d',\tilde a)\in u$.  The proof is similar, using Part (2b) in the Lemma.
\item For $(d',\tilde a)\in u$, $d'=0$ if and only if $\tilde a=0$.
  Namely, let $c\in C$ and $b''\in B''$
  such that $(d',c)\in h'\circ{\gamma'}^{-1}$, $(c,b'')\in g\circ\beta^{-1}$ and $(b'',\tilde a)\in f''\circ{\alpha''}^{-1}$.
  The following statements are equivalent:  $d'=0$; $d'\in\Im h'\circ 0$;
  $c\in\Im\gamma'g'=\Im g\beta'$ (use (1) in the Lemma, center
  at $C'$); $b''\in\Im\beta f=\Im f''\alpha$
  (again by (1), center at $B$); $\tilde a\in\Im \alpha''\circ0$; $\tilde a=0$.
\end{enumerate}

\smallskip
Again, (1) and (3) define a map $\Ker \tilde f\to \Ker\delta'$;
this map is onto by (2) and one-to-one by (3).
\end{proof}

\bigskip
The Proposition applies in particular to commutative diagrams given by
``multiplying'' a short right exact sequence and a short left exact sequence using the $\Hom$-functor.

\begin{theorem}\mylabel{theorem-isomorphic-homology}
Suppose $\mathcal E:X\stackrel u\to Y\to Z\to 0$ is a short right exact sequence
and $\mathcal A:0\to A\to B\stackrel g\to C$ a short left exact sequence.
The complexes
$$0\to \Cok(u,A)\to \Cok(u,B)\to \Cok(u,C)\to 0$$
and
$$0\to \Cok(Z,g)\to \Cok(Y,g)\to \Cok(X,g)\to 0$$
have isomorphic homology in each position.
\end{theorem}

\begin{proof}
The following diagram is commutative with exact rows and columns.
$$
\xymatrixcolsep{5mm}
\xymatrix{  & 0 \ar[d] & 0 \ar[d] & 0 \ar[d] &  &  \\
  0 \ar[r] & (Z,A) \ar[r] \ar[d] & (Y,A) \ar[r]^{(u,A)} \ar[d] & (X,A) \ar[r] \ar[d] & \Cok(u,A) \ar[r] & 0 \\
  0 \ar[r] & (Z,B) \ar[r] \ar[d]^{(Z,g)} & (Y,B) \ar[r]^{(u,B)} \ar[d]^{(Y,g)} & (X,B) \ar[r] \ar[d]^{(X,g)} & \Cok(u,B) \ar[r] & 0 \\
  0 \ar[r] & (Z,C) \ar[r] \ar[d] & (Y,C) \ar[r]^{(u,C)} \ar[d] & (X,C) \ar[d] \ar[r] & \Cok(u,C) \ar[r] & 0 \\
  & \Cok(Z,g) \ar[d] & \Cok(Y,g) \ar[d] & \Cok(X,g) \ar[d] & & \\
  & 0 & 0 & 0 & & }
$$

The result follows from Proposition~\ref{proposition-isomorphic-homology}.
\end{proof}

\section{Hammock Functions}
\label{sec-hammock}

In this section we give the proofs for Theorem~\ref{theorem-isolated} and Proposition~\ref{prop-hammock-factoring}
from the Introduction.  Then we state the dual results.

\begin{proof}[Proof of Theorem~\ref{theorem-isolated}]
Let $B\stackrel g\to C$ be a sink map for $C$ in $\mathcal C$, so
the expression $\dim\Cok(M,g)$ determines the multiplicity of $C$ as a 
summand of $M\in\mathcal C$.  (Namely, if $d=\dim_k\End(C)/\rad\End(C)$ then
$\dim_k\Cok(M,g)=d\cdot \mu_C(M)$.)
If $C$ is an indecomposable non-projective object in $\mathcal C$, then it occurs
as the end term of an Auslander-Reiten sequence 
$\mathcal A: 0\to A\to B\stackrel g\to C\to 0$ in $\mathcal C$.
In case $C$ is indecomposable projective in $\mathcal C$, the sink map $g$
is still part of a short left exact sequence $\mathcal A: 0\to A\to B\stackrel g\to C$.

\smallskip
In the statement of the theorem,
the short right exact sequence is given by $\mathcal E:X\stackrel u\to Y\to Z\to 0$.
We apply Theorem~\ref{theorem-isomorphic-homology} to $\mathcal A$ and $\mathcal E$,
then the following two complexes have isomorphic homology.
$$
\xymatrixrowsep{.2cm}
\xymatrix{0 \ar[r] & \Cok(u,A) \ar[r] & \Cok(u,B) \ar[r] & \Cok(u,C) \ar[r] & 0 \\
  0 \ar[r] & \Cok(Z,g) \ar[r] & \Cok(Y,g) \ar[r] & \Cok(X,g) \ar[r] & 0 }
$$
Note that the dimensions of the vector spaces in the first sequence are the values
of the hammock function.  Up to a factor of $d$,
the dimensions of the modules in the second sequence are the multiplicities
of $C$ as a direct summand of $Z$, $Y$, or $X$, respectively.

\smallskip
1. The hammock function $e_*$ has an isolated source $C$ by definition 
if $\dim\Cok(u,C)>\dim\Cok(u,B)$.
In this case, $\Cok(X,g)\neq 0$, so $X$ has a summand isomorphic to $C$.

\smallskip
2. Assume that $A$ is not injective.
Then $e_*$ has an isolated sink $A$ if $\dim\Cok(u,A)>\dim\Cok(u,B)$.
In this situation, $Z$ has a summand isomorphic to $C=\tau^{-1}A$.

\smallskip
3. Moreover, $e_*$ has a tangent point at a mesh $\mathcal M: A\to \bigoplus B_i\to C$
where $C$ is indecomposable
if $\sum\dim\Cok(u,B_i)>\dim\Cok(u,A)+\dim\Cok(u,C)$.
Hence $Y$ has a summand isomorphic to $C$, the end point of the mesh.

\smallskip
4. The condition that $C$ does not occur as a direct summand of $X\oplus Y\oplus Z$,
is equivalent to the second sequence being constant zero, which means that 
the first sequence is exact, so $e_*$ is additive on $\mathcal A$.
\end{proof}

\medskip
\begin{proof}[Proof of Proposition~\ref{prop-hammock-factoring}]
  Let $\mathcal E:0\to X\stackrel u\to Y\to Z\to 0$ be the short non-split exact sequence
  in the Proposition, and let $f:X\to M$ be a morphism which
  does not factor over $u$.  Since $\mathcal E$ is non-split,
  the indecomposable module $Z$ is not a projective object in $\mathcal C$
  and the Auslander-Reiten sequence $\mathcal A:0\to \tau Z\to W\to Z\to 0$ in $\mathcal C$ is defined.
  We show that there exists a morphism $g:M\to \tau Z$
  such that $gf:X\to \tau Z$ does not factor over $u$.
  $$
  \xymatrix{\mathcal A: & 0 \ar[r] & \tau Z \ar[r]^{w'} & W \ar[r]^w & Z \ar[r] & 0 \\
    \mathcal E: & 0 \ar[r] & X \ar[u]^h \ar[r]^u \ar[d]_f & Y \ar[u]^{h'} \ar[r]^v \ar[d]_{f'} \ar[ul]_d & Z \ar@{=}[u] \ar[r] \ar@{=}[d] & 0 \\
    \mathcal Ef: & 0 \ar[r] & M \ar[r]^{p'} \ar@/_1pc/[uu]_<<<g & P \ar[r]^p \ar@/_1pc/[uu]_<<<{g'} & Z \ar[r] & 0 }
  $$
  Since $v:Y\to Z$ is not a split epimorphism, there exists $h':Y\to W$ such that
  $v=wh'$.  Then the kernel map $h:X\to \tau Z$ does not factor over $u$
  since otherwise $1_Z$ would factor over $w$ which is not possible since $\mathcal A$
  is non-split.  In particular, the target $\Cok(u,\tau Z)$ of the map in the
  Proposition is non-zero.

  \smallskip
  Since $f$ does not factor over $u$, the map $1_Z$ does not factor over $p$
  so the induced sequence $\mathcal Ef$ is non-split.
  Since $\mathcal C$ is extension-closed, the middle term in $\mathcal Ef$ is
  an object in $\mathcal C$. 
  Hence there is a map $g':P\to W$ with $p=wg'$. Let $g:M\to\tau Z$ be the kernel map.
  We show $gf-h$ factors over $u$, so $gf$ does not factor over $u$.

  \smallskip
  Since $w (g'f'-h') = pf'-v= v-v = 0$, there is a map $d:Y\to \tau Z$  with
  $w'd=g'f'-h'$.  Then $w' (gf-h)=g'p'f-h'u=(g'f'-h')u=w'du$ implies that
  $gf-h=du$ since $w'$ is a monomorphism.  In particular, $gf-h$ factors over $u$.
  We have seen that $h$ does not factor over $u$, so $gf$ does not factor over $u$.
\end{proof}

\bigskip
We state the dual results for Theorem~\ref{theorem-isolated} and Proposition~\ref{prop-hammock-factoring}.

\medskip
Suppose  $E:\mod\Lambda\to \mod k$ is a finitely presented contravariant functor,
say with projective resolution
$$0\longrightarrow (-,X)\longrightarrow(-,Y)\longrightarrow(-,Z)\longrightarrow E\longrightarrow 0,$$
given by a short left exact sequence $\mathcal E:0\to X\to Y\stackrel v\to Z$.
We call the map
$$e^*:\Gamma_0\to \mathbb N_0, \;A\mapsto \dim\Hom\Cok(A,v)$$
the \boldit{hammock function} for the contravariant functor $E$ (or for the sequence $\mathcal E$).

\medskip
Here is the dual version of Theorem~\ref{theorem-isolated}.

\begin{theorem}
  \mylabel{theorem-isolated-dual}
  Let $\mathcal E:0\to X\to Y\stackrel v\to Z$ be a short left exact sequence and
  $e^*:\Gamma_0\to\mathbb N_0$ the corresponding hammock function.
  \begin{enumerate}
  \item If $A$ is an isolated sink for $e^*$ then $A$ is a direct summand of $Z$.
  \item If $C$ is an isolated source for $e^*$ then $C$ is projective in $\mathcal C$
    or $A=\tau C$ is a direct summand of $X$.
  \item If the mesh $\mathcal M:A\to \bigoplus B_i\to C$ is an isolated tangent for $e^*$
    then $C$ is indecomposable projective or $A$ is an indecomposable direct summand of $Y$.
  \item If $\mathcal M:A\to \bigoplus B_i\to C$ is a mesh such that $A$ is
    indecomposable and does not occur as a direct summand of $X\oplus Y\oplus Z$
    then the hammock function $e^*$ is additive on $\mathcal M$.\qed
  \end{enumerate}
  
\end{theorem}

For the dual version of Proposition~\ref{prop-hammock-factoring},
suppose $\mathcal E:0\to X\to Y\stackrel v\to Z\to 0$ is a non-split short exact sequence
with $X$ indecomposable non-injective.  Then the modules $M$ in the hammock
are characterized by lying between $\tau^{-1}X$ and one of the summands of $Z$ in the sense
that there is a map $f:M\to Z$ which does not factor over $v$ and a map $g:\tau^{-1}X\to M$
such that the product $fg:\tau^{-1}X\to Z$ does not factor over $v$:

\begin{proposition}
  \mylabel{prop-hammock-factoring-dual}
  Let $\mathcal E:0\to X\to Y\stackrel v\to Z\to 0$ be a non-split short exact sequence
  in $\mathcal C$ with
  $X$ indecomposable.  Then for any $M\in\mathcal C$, the bilinear form given by composition,
  $$\Cok(M,v)\times \Hom(\tau^{-1}X,M) \longrightarrow \Cok(\tau^{-1}X, v),$$
  is left non-degenerate.\qed
\end{proposition}

\begin{example}
  We illustrate the two results in the introduction and their dual versions
  in the example of a hammock that contain exactly one vertex.
  Let $M$ be an indecomposable module, neither projective nor injective in $\mathcal C$.

  \smallskip
  For Theorem~\ref{theorem-isolated}, consider the Auslander-Reiten sequence 
  $\mathcal E:0\to X\stackrel u\to Y\to Z\to 0$ starting at $X=M$.
  The hammock function $e_*$ counts the multiplicity of $M$ as a direct summand of its object,
  hence $e_*$ vanishes on $\Gamma$, except at the vertex $M$.
  In fact, $M=X$ is an isolated source; for each indecomposable summand $D$ of $Y$, the
  mesh ending at $D$ (which contains $M$ as a summand of the middle term) is an isolated tangent;
  and the sink of the hammock is also $M$ since $M=\tau Z$.

  \smallskip
  Regarding Proposition~\ref{prop-hammock-factoring}, note that $M$ is the only object in $\Gamma$
  between $X$ and $\tau Z$.

  \smallskip
  For Theorem~\ref{theorem-isolated-dual}, take the Auslander-Reiten sequence
  $\widetilde{\mathcal E}: 0\to \tilde X\to \tilde Y\stackrel v\to \tilde Z\to 0$ ending at $\tilde Z=M$.
  Then $e^*$ vanishes on $\Gamma$ except at $M$.
  We have $M$ as an isolated sink since $M=\tilde Z$;
  for each indecomposable summand $D$ of $\tilde Y$, the mesh starting at $D$ (which contains $M$ as a
  summand of the middle term) is an isolated tangent; and the source of the hammock
  is $M$ since $M=\tau^{-1}\tilde X$.

  \smallskip
  This also illustrates Proposition~\ref{prop-hammock-factoring-dual} since $M$ is the only object
  in $\Gamma$ between $\tau^{-1}\tilde X$ and $\tilde Z$.
\end{example}

\section{Applications I:  Modules}
\label{sec-quiver}

In this section we use hammocks to visualize some properties of modules,
in particular properties related to multiplication functions.

\subsection{The multiplication map for modules}

Let $R$ be an algebra, $M$ a finite dimensional (right) $R$-module and $a\in R$.
The multiplication map
$$\mu_{M,a}:\; M\to M, \; m\mapsto ma,$$
gives rise to subspaces of $M$, in particular, $\Ker\mu_{M,a}$, $\Im\mu_{M,a}$ and $\Cok\mu_{M,a}$,
the dimension of which is a hammock function.

\medskip
As a left $\End M$-module, the \boldit{kernel} of the multiplication map is given as a Hom space by the isomorphism
$$\Hom_R(R/aR,M)\;\cong\;\Ker\mu_{M,a},\quad f\mapsto f(1).$$
Hence the dimension of the kernel is the covariant defect $\dim\Ker\mu_{M,a}=e_*(M)$ where
$\mathcal E$ is the short right exact sequence
$$\mathcal E: \quad R/aR\longrightarrow 0\longrightarrow 0\longrightarrow 0.$$
Considering Theorem~\ref{theorem-isolated}, the hammock function has sources at the indecomposable
direct summands of $R/aR$, no non-injective isolated sinks, and no tangents with indecomposable end point.

\smallskip
As a special case, consider $a=1-e$ for an idempotent $e\in R$.  Then the above hammock function
yields $\dim Me=\dim\Ker\mu_{M,a}$. 
The arising hammocks have already been studied in \cite{b}.

\medskip
The \boldit{image,} being isomorphic to the factor of $M$ modulo the kernel of the multiplication map,
is obtained as the cokernel of the sequence
$$
\xymatrix{0 \ar[r] & (R/aR,M) \ar[r] \ar[d]^\cong & (R,M) \ar[r] \ar[d]^\cong & E(M) \ar[r] \ar[d]^\cong & 0 \\
 0 \ar[r]  & \Ker\mu_\alpha \ar[r]^{{\rm incl}} & M\ar[r] & M/\Ker\mu_\alpha\cong\Im\mu_\alpha \ar[r] & 0 }
$$
(where we abbreviate $\Hom(A,B)$ by $(A,B)$), 
thus $\dim\Im\mu_{M,a}=e_*(M)$ is the covariant defect of the short right exact sequence
$$\mathcal E:\quad  R\stackrel u\longrightarrow R/aR \longrightarrow 0\longrightarrow 0.$$

\medskip
For the \boldit{cokernel} of the multiplication map, consider the short right exact sequence
$$\mathcal E:\quad R\stackrel{u}\longrightarrow R\longrightarrow R/aR\longrightarrow 0$$
where $u$ is the left multiplication by $a$, which gives rise to the commutative diagram
$$
\xymatrix{ 0 \ar[r] & (R/aR,M) \ar[r] \ar[d]^\cong & (R,M) \ar[r] \ar[d]^\cong & (R,M) \ar[r] \ar[d]^\cong & E(M) \ar[r] \ar[d]^\cong & 0 \\
 0\ar[r]  & \Ker\mu_\alpha \ar[r]^{{\rm incl}} & M \ar[r]^{\mu_\alpha} & M \ar[r] & \Cok\mu_\alpha \ar[r] & 0 }
$$
hence $\dim\Cok\mu_{M,a}=e_*(M)$, as desired.

\subsection{Quiver representations}
\mylabel{sec-cok}
We use hammocks to visualize some data derived from quiver representations.

\medskip
Let $\Lambda=kQ$ be the path algebra
of a quiver $Q=(Q_0,Q_1,s,t)$.  For a $k$-linear representation $M$
of $Q$, a vertex $v\in Q_0$, an arrow $\alpha\in Q_1$, and a path $\pi=\alpha_1\cdots\alpha_s$ in $Q$,
we denote by $M_v$ the
vector space in position $v$, by $M_\alpha$ the linear map $M_{s(\alpha)}\to M_{t(\alpha)}$,
and by $M_\pi$ the linear map $M_{\alpha_s}\cdots M_{\alpha_1}:M_{s(\pi)}\to M_{t(\pi)}$.

\medskip
In particular, we consider the path algebra of the quiver $\mathbb E_6$ in bipartite orientation
$$Q: \qquad 
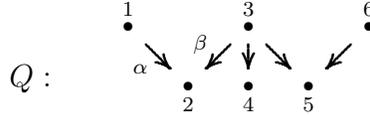
\beginpicture\setcoordinatesystem units <.8cm,.8cm>
\multiput{$\scriptstyle \bullet$} at 0 1  1 0  2 1  2 0  3 0  4 1 /
\arr{0.3 0.7}{0.7 0.3}
\arr{1.7 0.7}{1.3 0.3}
\arr{2 0.7}{2 0.3}
\arr{2.3 0.7}{2.7 0.3}
\arr{3.7 0.7}{3.3 0.3}
\put{$\scriptstyle 1$} at 0 1.3
\put{$\scriptstyle 2$} at 1 -.3
\put{$\scriptstyle 3$} at 2 1.3
\put{$\scriptstyle 4$} at 2 -.3
\put{$\scriptstyle 5$} at 3 -.3
\put{$\scriptstyle 6$} at 4 1.3
\put{$\scriptstyle \alpha$} at 0.2 0.3 
\put{$\scriptstyle \beta$} at 1.2 0.7
\endpicture$$
In this section, we picture the hammocks given by the kernel of the linear map
$M_\beta:M_3\to M_2$, by the image of $M_\beta$, and by the cokernel of $M_\alpha$.

\medskip
First, we consider \boldit{the kernel at an arrow or a path.}
For vertices $i,j\in Q_0$ and an arrow $\alpha:i\to j$ in $Q_1$, denote by
$P(i)$ the indecomposable projective module corresponding to vertex $i$
and by $(\alpha)$ the submodule of $P(i)$ given by
the image of the map $P(\alpha):P(j)\to P(i)$; similarly, for a path $\pi=\alpha_1\cdots\alpha_s:i\to j$,
we write $(\pi)$ for the submodule of $P(i)$ given by the image of the composition
$P(\pi)=P(\alpha_1)\cdots P(\alpha_s):P(j)\to P(i)$.  We claim that the short right exact sequence
$$\mathcal E: \quad P(i)/(\pi)\stackrel v\longrightarrow 0\longrightarrow 0\longrightarrow 0$$
defines the hammock function $e_*(M)=\dim\Ker M_\pi$.

\smallskip 
Namely, applying the functor $\Hom(-,M)$ to the short exact sequence
  $0\to P(j)\to P(i)\to P(i)/(\pi)\to 0$ gives the long exact sequence starting with
  $$0\to (P(i)/(\pi),M)\to (P(i),M) \stackrel g\to (P(j),M)\to \cdots.$$
  Identifying $(P(i),M)$ with $M_i$, $(P(j),M)$ with $M_j$ and $g$ with $M_\pi$,
  we obtain $\Ker M_\pi\cong \Hom(P(i)/(\pi),M)$ and the claim follows.

  \smallskip
  We obtain from Theorem~\ref{theorem-isolated}:
  In the Auslander-Reiten quiver, $e_*$ is the hammock function with
  source $P(i)/(\pi)$, no tangents with indecomposable end point, and no non-injective sinks.

\begin{example}
  Figure~\ref{fig-ker-beta} shows the Auslander-Reiten quiver for the path algebra $kQ$ of type $\mathbb E_6$
  given above.  Each object is represented by its dimension vector.
  The region encircled by the dotted curve is the hammock for
  the function $e_*(M)=\dim \Ker M_\beta$.
  The hammock has source $X=P(3)/(\beta)$, no tangents containing the end point, but there
  is a tangent at the mesh starting at $\tilde Y$, and no non-injective sinks, but there is
  an injective sink at $\tilde Z=I(3)$.

\begin{figure}[ht]
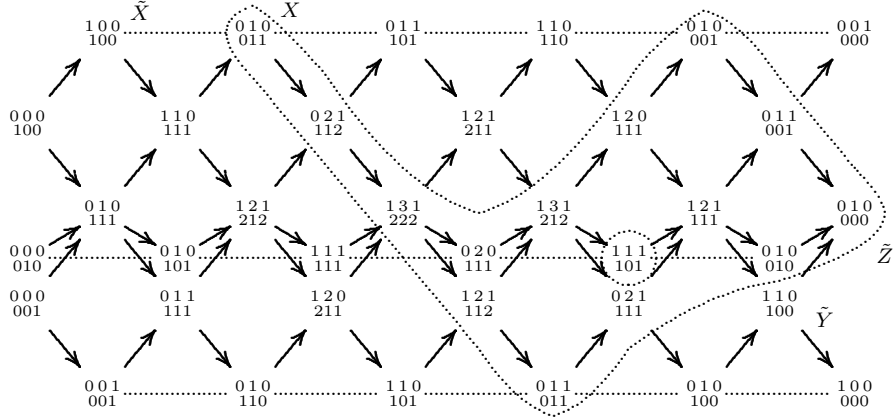

$$
\beginpicture\setcoordinatesystem units <1cm,.6cm>
\arqesix
\put{$\scriptstyle X$} at 3.5 8.5
\put{$\scriptstyle \tilde X$} at 1.5 8.5
\put{$\scriptstyle \tilde Y$} at 10.6 1.7
\put{$\scriptstyle \tilde Z$} at 11.4 3.1
\setquadratic
\setdots<2pt>
\plot 4 5  5 3  6 1  6.5 0  7 -.5  7.5 0  8 1  9 2  10 2.5  10.5 2.8  11 3.2  11.3 3.6
11.4 4  11.3 4.4  11 5  10.5 6  10 7  9.5 8  9 8.5  8.5 7.9  8.1 7.2  7 5  6 4  5 5  4 7
3.45 8  3 8.6  2.65 8  3 7  3.5 6  4 5 /
\circulararc 360 degrees from 8 2.4 center at 8 3
\endpicture
$$
\caption{The hammock for the functor $E:M\mapsto \Ker M_\beta$}
\label{fig-ker-beta}
\end{figure}

  The path $\pi$ also defines a canonical map between indecomposable injective modules, $v:I(j)\to I(i)$.
  Using the contravariant defect and Theorem~\ref{theorem-isolated-dual}, one can obtain
  $\Ker M_\pi\cong \Cok(M_\pi^*)\cong \Cok\Hom(M,v)$.
  In the above example, take $\pi=\beta$ and
  consider the short left exact sequence
  $$\widetilde{\mathcal E}:0\to P(2)\to I(2)\stackrel v\to I(3).$$
  Then $\tilde Z=I(3)$ is the sink for the hammock function $e^*$; there is a tangent
  at the mesh starting at $\tilde Y=I(2)$; and the source occurs at $\tau^{-1}\tilde X$
  where $\tilde X=\Ker v\cong P(2)$.  The modules $\tilde X$, $\tilde Y$, $\tilde Z$ are
  indicated in the above Auslander-Reiten quiver.
\end{example}

\medskip
We consider \boldit{the image at an arrow or at a path.}
For $\alpha:i\to j$ an arrow (or $\pi=\alpha_1\cdots\alpha_s:i\to j$ a path),
consider the short right exact sequence
$$\mathcal E:\quad P(i)\stackrel{u}\longrightarrow P(i)/(\alpha) \longrightarrow 0
\longrightarrow 0,$$ and the exact sequence obtained by applying the functor
$\Hom(-,M)$.
$$
\xymatrix{ 0 \ar[r] & (P(i)/(\alpha),M) \ar[r]^{(u,M)} \ar[d]^\cong & (P(i),M) \ar[r] \ar[d]^\cong & \Cok(P(\alpha),M) \ar[r] \ar[d]^\cong & 0\\
 0 \ar[r]  &\Ker M_\alpha \ar[r]^{{\rm incl}} & M_i \ar[r]^{{\rm can}} & \Im M_\alpha \ar[r] & 0}
$$

Thus, the dimension of the image of the map associated with an arrow $\alpha$
is given by a hammock function with source $P(i)$, tangent at $P(i)/(\alpha)$
(or, in case of a path $\pi$, at $P(i)/(\pi)$), and no non-injective sink.

\begin{figure}[ht]
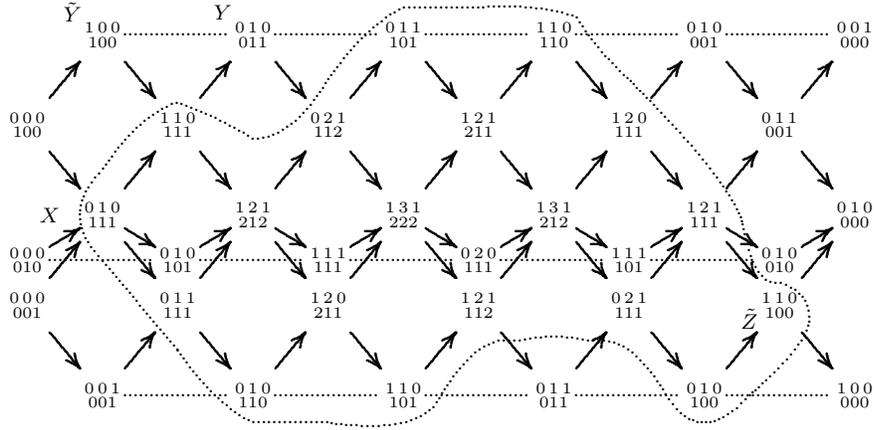

$$
\beginpicture\setcoordinatesystem units <1cm,.6cm>
\arqesix
\put{$\scriptstyle X$} at 0.3 4
\put{$\scriptstyle Y$} at 2.6 8.5
\put{$\scriptstyle \tilde Y$} at .6 8.5
\put{$\scriptstyle \tilde Z$} at 9.6 1.7
\setquadratic
\setdots<2pt>
\plot 2 1  2.5 0  3 -.6  3.5 -.6  4 -.6  5 -.6  5.5 0  6 1  6.5 1.2  7 1.3
7.5 1.2  8 1  8.5 0  9 -.6  9.6 0  10.2 1  10.4 1.7  10.3 2.2  10 2.5  9.7 2.8  9.5 3.9
9 5.2  8 7  7.6 7.7  7.3 8.4  6.7 8.6  6 8.6  5.5 8.6  5 8.6  4.5 8  4 7  3.5 6  3 5.7  2.5 6.15
2 6.5  1.5 6  1 5  0.7 4  1 3  1.5 2  2 1 /
\endpicture
$$
\caption{The hammock for the functor $E: M\mapsto\Im M_\beta$}
\label{figure-im-beta}
\end{figure}

\begin{example}
  The indecomposable modules $M$ in the Auslander-Reiten quiver for $kQ$
  for which $\Im M_\beta\neq 0$ are pictured in Figure~\ref{figure-im-beta}.
  The hammock is given by the
  sequence $\mathcal E:\; P(3)\stackrel{u}\longrightarrow P(3)/(\beta)\to 0\to 0$,
  it has source $X=P(3)$, one tangent before $Y=P(3)/(\beta)$
  and injective sink $I(2)$.
  Note that there are five indecomposables $M$ in the hammock with
  $e_*(M)=\dim\Im M_\beta=2$; they are characterized here by $\dim M_2=2$.
\end{example}

The same hammock can be obtained dually using Theorem~\ref{theorem-isolated-dual}
and the short left exact sequence
$\widetilde{\mathcal E}: 0\to 0\to \Ker_{I(2)}(\beta)\stackrel{\rm incl}\to I(2)$.
Thus, the hammock has a sink at $\tilde Z=I(2)$, a tangent at the mesh starting
at $\tilde Y=\Ker_{I(2)}(\beta)$, and a source at some projective module or modules.

\medskip
Next, \boldit{the cokernel at an arrow or at a path.}
For $\alpha:i\to j$ an arrow, the right exact sequence
$$\mathcal E: \quad P(j)\stackrel v\longrightarrow P(i)\longrightarrow
P(i)/(\alpha)\longrightarrow 0$$
where $v=P(\alpha)$ gives rise to the top sequence in the diagram
$$
\xymatrix{\Hom(P(i),M) \ar[r]^{(v,M)} \ar[d]^\cong & \Hom(P(j),M)\ar[r]\ar[d]^\cong & \Cok(v,M)\ar[r] \ar[d]^\cong & 0 \\
  M_i \ar[r]^{M_\alpha} & M_j \ar[r] & \Cok M_\alpha \ar[r] & 0 }
$$
Thus, $\dim\Cok M_\alpha=e_*(M)$ is given by the hammock function for $\mathcal E$.
The hammock has a source at $P(j)$, a tangent at the mesh ending at $P(i)$
and a sink at $\tau [P(i)/(\alpha)]$.

\begin{example}
  \sloppy
  In the Auslander-Reiten quiver in Example~\ref{ex-intro}, the hammock consists
  of those indecomposable $kQ$-modules $M$ for which $\Cok M_\alpha\neq 0$.
  It is a hammock given by the sequence
  $\mathcal E:\; P(2)\stackrel v\to P(1)\to P(1)/(\alpha)\to 0.$
  The hammock has source at $X=P(2)$, one tangent before $Y=P(1)$ and a
  non-injective sink at $Z=\tau P(1)/(\alpha)=\tau I(1)$, as indicated.  
\end{example}

\subsection{The homology at two composable arrows or paths}

Let $\pi:i\to j$ and $\rho:j\to k$ be paths in $Q$, and $M$ a representation of $Q$
such that the composition $M_{\pi\rho}$ of $M_\pi:M_i\to M_j$ and $M_\rho:M_j\to M_k$ is the zero map.

\medskip
We are interested in the subquotient $\Ker M_\rho/\Im M_\pi$ of $M_j$.

\medskip
The following sequence is short right exact.
$$\mathcal E:\quad P(j)/(\rho) \stackrel{P(\pi)}\longrightarrow P(i)/(\pi\rho)\longrightarrow P(i)/(\pi)
\longrightarrow 0$$

Applying the functor $\Hom(-,M)$ we obtain
\begin{eqnarray*}
0\to (P(i)/(\pi),M) 
\to (P(i)/(\pi\rho),M) & \stackrel{(P(\pi),M)}\longrightarrow & (P(j)/(\rho),M)\\
                       & \to                                   &  \Cok (P(\pi),M)\to0
\end{eqnarray*}
where the first and second modules are isomorphic to $\Ker M_\pi$ and $M_i=\Ker M_{\pi\rho}$;
and the third to $\Ker M_\rho$.  Thus,
\begin{eqnarray*}
  \dim \Cok(P(\pi),M) & = & \dim \Ker M_\rho-(\dim M_i-\dim \Ker M_\pi)\\
  & = & \dim \Ker M_\rho/\Im M_\pi.
\end{eqnarray*}

Thus, the hammock given by the homology of two composable paths $\pi$ and $\rho$ has source at $P(j)/(\rho)$,
a unique non-injective sink at $\tau P(i)/(\pi)$ and a unique tangent with indecomposable end term at
$P(i)/(\pi\rho)$.

\section{Applications II:  Invariant subspaces}
\mylabel{sec-invariant}

In this section, let $R=k[T]/(T^n)$ be the bounded polynomial ring,
then the finite dimensional $R$-modules form the category $\mathcal N(n)=\mod R$
of nilpotent linear operators with nilpotency index at most $n$.
Each indecomposable $R$-module has the form $P^\ell=k[T]/(T^\ell)$ for some $1\leq \ell\leq n$.

\medskip
The category $\mathcal S(n)$ or $\mathcal S(R)$ of invariant subspaces of nilpotent linear operators
has as objects the pairs $(V,U)$ where $V\in\mathcal N(n)$
and $U$ is a submodule of $V$.
For example, if $\ell\leq m$, then the pair $P^m_\ell=(P^m,\soc^\ell P^m)$ is an object
in $\mathcal S(n)$. 
Such objects, called pickets, give rise to many interesting hammock functions.
They also turn out to be decisive in the study of the geometry of the representation space,
see \cite{ks}.

\medskip
We picture an object $V=e_1P^{\ell_1}\oplus\cdots e_sP^{\ell_s}$ in $\mathcal N(n)$ as $s$ columns
of boxes, where the $i$-th column has length $\ell_i$.  Thus, the $j$-th box in the $i$-th column
represents the basic element $e_iT^{j-1}$ in $V$.
Each submodule $U$ of $V$ considered in this paper has the property that $U$ has a set of generators
which are sums (not arbitrary linear combinations) of basic elements of $V$.
Hence we can picture each generator as a connected row of dots.
Here are some examples.

\newcommand\Pfour{\begin{picture}(3,12)
    \multiput(0,0)(0,3)4{\smbox}
  \end{picture}
}
\newcommand\Pthreefour{\begin{picture}(3,12)
    \multiput(0,0)(0,3)4{\smbox}
    \put(1.5,7.5){\sbullet}
  \end{picture}
}
\newcommand\Pbig{\begin{picture}(6,12)
    \multiput(0,0)(0,3)4{\smbox}
    \multiput(3,3)(0,3)2{\smbox}
    \put(1.5,7.5)\sbullet
    \put(4.5,7.5)\sbullet
    \put(4.5,4.5)\sbullet
    \put(1.5,7.5){\line(1,0)3}
  \end{picture}
}

$$P^4:\raisebox{-5mm}{\Pfour},\quad
P_3^4:\raisebox{-5mm}{\Pthreefour}, \quad
P_{31}^{42}=(e_1P^4\oplus e_2P^2,\langle e_1T+e_2,e_2T\rangle):\raisebox{-5mm}{\Pbig}.$$

\medskip
The triangular matrix ring $\Lambda=\big(\begin{smallmatrix} R & R\\ 0 & R \end{smallmatrix}\big)$
is a Gorenstein algebra, so the Gorenstein-projective objects form a full exact extension-closed
subcategory which has Auslander-Reiten sequences.  This category is equivalent to $\mathcal S(n)$,
see \cite{lz}.  Thus we are dealing with a Frobenius category with two indecomposable projective-injective
objects, $P_0^n$ and $P_n^n$.

\medskip
In Figure~\ref{figure-gamma-s4},
we picture the Auslander-Reiten quiver $\Gamma$ for $\mathcal S(4)$, taken from \cite[(6,4)]{rs}.
There are 20 indecomposable objects in $\Gamma$;
note that the quiver is periodic.  We picture here one fundamental domain, so 
the objects on the left and on the right are to be identified, up to a reflection.

\def\scale#1{\makebox[0pt][c]{\beginpicture\setcoordinatesystem units <1.8mm,1.8mm>#1
    \endpicture}}
\def\arqfouronepic{\beginpicture\setcoordinatesystem units <.55cm,.45cm>
  \multiput{} at -4 0  13 12 /
  \arqfourone
  \endpicture
}
\begin{figure}[ht]
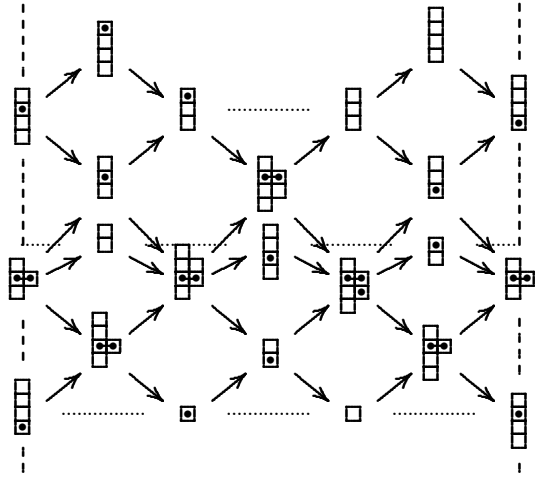

  $$  \arqfouronepic$$
  \caption{The Auslander-Reiten quiver for $\mathcal S(4)$}
  \label{figure-gamma-s4}
\end{figure}

\medskip
In this section, we consider hammock functions which, for an object $(V,U)\in\mathcal S(n)$,
count the multiplicity of a module $P^\ell\in\mathcal N(n)$ as a direct summand of $V$
(or of $U$, or of $V/U$).  In our last example, the hammock function counts
the number of indecomposable direct summands of $U$.
  
\subsection{Our tool:  Adjoint functors}

In this section we assume that $\mathcal C$ and $\mathcal D$ are exact $k$-categories with Auslander-Reiten
sequences.

\begin{proposition}
  \label{prop-RA}
  Suppose $R$ and $L$ form  a pair of adjoint functors
  $$
  \xymatrix{ \mathcal C\quad \ar@/^/[r]^R & \quad\mathcal D \ar@/^/[l]^L }
  $$
  and $M$ is an indecomposable object in $\mathcal D$ with source map $s:M\to N$.
  Put $d=\dim\End M/\rad\End M$.  The short right exact sequence
  $$\mathcal E:\quad LM\stackrel{Ls}\longrightarrow LN\longrightarrow L\Cok(s)\longrightarrow 0$$
  gives rise to a hammock function $e_*$ such that, for an object $A\in \mathcal C$,
  the multiplicity of $M$ as a direct summand of $RA$ is $\frac 1d e_*(A)$
\end{proposition}

\begin{proof}
  The source map $s:M\to N$ gives rise to the short right exact sequence in $\mathcal D$,
  $$M\stackrel s\longrightarrow N\longrightarrow \Cok(s)\longrightarrow 0$$
  and hence, since a left adjoint functor is always right exact, to the short right exact
  sequence in $\mathcal C$,
  $$\mathcal E: \quad LM\stackrel{Ls}\longrightarrow LN\longrightarrow L\Cok(s)\longrightarrow 0.$$
  We can now compute using the adjoint isomorphism
  \begin{eqnarray*}
    d\cdot\mu_M(RA) & = & \dim\Cok\Hom_{\mathcal D}(s,RA)\\
    & = & \hom_{\mathcal D}(\Cok(s),RA)-\hom_{\mathcal D}(N,RA)+\hom_{\mathcal D}(M,RA)\\
    &= & \hom_{\mathcal C}(L\Cok(s),A)-\hom_{\mathcal C}(LN,A)+\hom_{\mathcal C}(LM,A)\\
    &=& \dim\Cok\Hom_{\mathcal C}(Ls,X).
  \end{eqnarray*}
  
\end{proof}

\subsection{The multiplicity of $P^m$ as a summand of $V$}
\mylabel{sec-ambient}

Consider the categories $\mathcal C=\mathcal S(n)$ and $\mathcal D=\mathcal N(n)$,
and the functors $R:\mathcal S(n)\to \mathcal N(n), (V,U)\mapsto V$ and
$L:\mathcal N(n)\to \mathcal S(n), V\mapsto (V,0)$ which form an adjoint pair:
$$\xymatrix{\mathcal S(4)\quad \ar@/^/[r]^{(U\subset V)\mapsto V} &
 \quad \mathcal N(4) \ar@/^/[l]^{(0\subset V)\gets V }}
$$

\medskip
Starting with an Auslander-Reiten sequence for $P^m$ in $\mathcal N(n)$,
$$(*)\qquad 0\to P^m\to P^{m+1}\oplus P^{m-1}\to P^m\to 0$$
if $m<n$ (where $P^0=0$ in case $m=1$) or with the source map $P^n\to P^{n-1}$,
we obtain by applying $L$ the corresponding sequence in $\mathcal S(n)$:
$$\mathcal E: \quad P^m_0\stackrel u\longrightarrow P^{m+1}_0\oplus P^{m-1}_0\longrightarrow P^m_0
\longrightarrow 0$$
for $m<n$ or $\mathcal E: P^n_0\stackrel u\to P_1^n$ if $m=n$.

\medskip
It follows from the Proposition that for an object $A=(V,U)$ in $\mathcal S(n)$,
the multiplicity of $P^m$ as a direct summand of $V$ is
$$\mu_{P^m}(V)\;=\; \dim\Cok\Hom(u,A)\;=\;e_*(A).$$

\begin{example}
  We picture in Figure~\ref{figure-summands}
  the hammock in the category $\mathcal S(4)$ given by those objects $(V,U)$ where
  $V$ contains a direct summand isomorphic to $P^3$.
  For clarity, we present two copies of the fundamental domain of the above Auslander-Reiten quiver.

  \smallskip
  The objects in the hammocks which we consider in sections \ref{sec-ambient}, \ref{sec-subspace} and \ref{sec-factor}
  are marked by a large semicircle on the left, it is solid, dotted or dashed, depending on the section:
  The ambient space is marked by a solid, the subspace by a dotted, and the factor by a dashed semicircle.

\def\scale#1{\makebox[0pt][c]{\beginpicture\setcoordinatesystem units <1.8mm,1.8mm>#1
    \endpicture}}
  \def\halfcirc{\circulararc 180 degrees from .2 .6  center at  .2 0 }  
  \def\smallhalfcirc{\circulararc 180 degrees from .2 .3  center at  .2 0 } 
\def\arqfoure{\beginpicture\setcoordinatesystem units <.45cm,.45cm>
  \multiput{} at -2 0  25 15 /
  \arqfourblack
  \setsolid
  \multiput{\begin{red}{$\halfcirc$}\end{red}} at 5.7 11  7.7 9  9.6 6  11.7 4  13.7 2 / %
  \multiput{\begin{red}{$\smallhalfcirc$}\end{red}} at 17.7 11  19.5 9  21.5 6  23.5 4  -.5 4  1.7 2 / %
  \setdots<2pt>
  \multiput{ \green{$\halfcirc$} } at 1.7 11  3.5 8.6  5.5 5.6  7.5 3.6  9.7 1.6 /
  \multiput{ \green{$\smallhalfcirc$} } at 13.7 11  15.7 8.6  17.5 5.6  19.7 3.6  21.7 1.6 /
  \setdashes<2pt>
  \multiput{ \blue{$\halfcirc$} } at 9.7 11.4  11.5 9.4  13.5 6.4  15.5 4  17.7 2 /
  \multiput{ \blue{$\smallhalfcirc$} } at 21.7 11.4  23.7 9.4 -.3 9.4  1.5 6.4  3.7 4  5.7 2 /
  \put{$\sssize X^s$} at 2 12
  \put{$\sssize Y_1^s$} at 8 6
  \put{$\sssize Y_2^s$} at 12.9 -.5
  \put{$\sssize Z^s$} at 14.5 1
  \put{$\sssize X^a$} at 6 12
  \put{$\sssize Y_2^a$} at 8.9 13.5
  \put{$\sssize Y_1^a$} at 12 6.3
  \put{$\sssize Z^a$} at 17.5 0.7
  \put{$\sssize X^f$} at 9.5 12.5
  \put{$\sssize Y_1^f$} at 16 5.5
  \put{$\sssize Y_2^f$} at 20.9 -.5
  \put{$\sssize Z^f$} at 22.5 1
  \endpicture
}

\begin{figure}[ht]
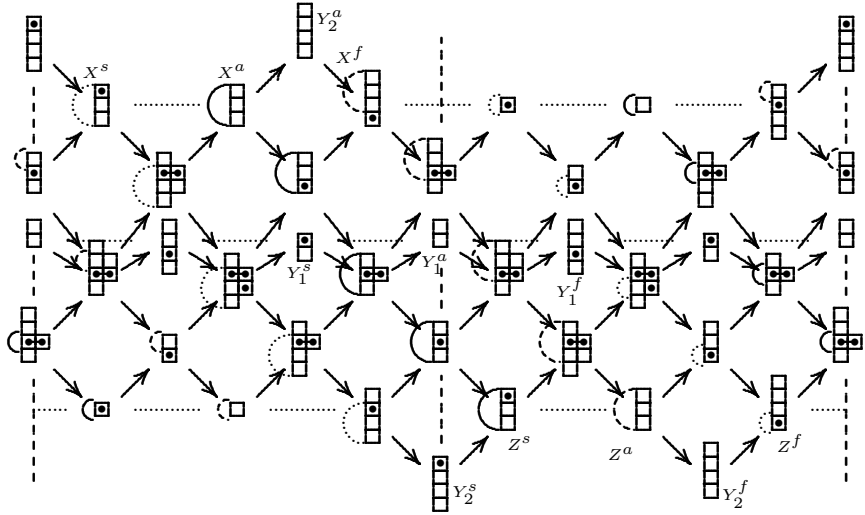

\begin{center}
  \arqfoure
\end{center}
\caption{Hammocks for summands of $U$, $V$ and $V/U$}
\label{figure-summands}
\end{figure}

The large semicircle indicates a summand of type $P^3$; some other objects are marked with a small semicircle,
which indicates a summand of type $P^1$.
If the semicircle is solid, the summand occurs in the ambient space $V$; if it is dotted, it occurs
in the subspace $U$; and the dashed semicircle represents a summand of the factor $V/U$.

\medskip
The five objects marked with a large solid semicircle on the left are characterized by having a summand
of type $P^3$ in the ambient space $V$.  The short right exact sequence
$\mathcal E:P_0^3\to P_0^2\oplus P_0^4\to P_0^3\to 0$ in this section
defines a hammock with source $X^a=P_0^3$, two tangents at the meshes ending at $Y^a_1=P_0^2$
and $Y_2^a=P_0^4$ 
and a sink at $\tau Z^a$ where $Z^a=P_0^3$.
There is another tangent at the mesh starting at the injective $P_4^4$ labeled $Y_2^s$ which has no indecomposable end term.

\smallskip
Note that in the picture, the semicircle indicates the position of $P^3$ as a direct summand
of the ambient space.
\end{example}

\subsection{The multiplicity of $P^m$ as a summand of $U$}
\mylabel{sec-subspace}

Here the functors  $R:\mathcal S(n)\to \mathcal N(n), (V,U)\mapsto U$,
and $L:\mathcal N(n)\to \mathcal S(n), U\mapsto (U,U)$ form an adjoint pair as follows:
$$\xymatrix{\mathcal S(4)\quad \ar@/^/[r]^{(U\subset V)\mapsto U} &
  \quad \mathcal N(4) \ar@/^/[l]^{(U = U)\gets U }}
$$

\medskip
The right exact sequence obtained from the above source map (*) in $\mathcal N(n)$ is
$$\mathcal E:\quad P^m_m\stackrel u\longrightarrow P^{m+1}_{m+1}\oplus P^{m-1}_{m-1}
\longrightarrow P^m_m\longrightarrow 0$$
if $m<n$ and $E: P^n_n\stackrel u\to P_{n-1}^{n-1}$ otherwise.

\medskip
We obtain from the Proposition that for an object $A=(V,U)$, the multiplicity
of $P^m$ as a direct summand of $U$ is
$$\mu_{P^m} (U)\;=\; \dim\Cok\Hom(u,A) \;=\; e_*(A).$$

\begin{example}
  Returning to the example where $n=4$ and $m=3$,
  the five objects in the Auslander-Reiten quiver marked with a large dotted semicircle
  have $P^3$ as a direct summand of the subspace.
  In fact, the dotted semicircle indicates the position of the subspace.
  The sequence in this section, $\mathcal E: P_3^3\to P_2^2\oplus P_4^4\to P_3^3\to 0$
  defines a hammock with source $X^s=P_3^3$, tangents at the meshes ending in $Y^s_1=P_2^2$ and $Y^s_2=P_4^4$;
  and a sink at $\tau Z^s$ where $Z^s=P_3^3$.  There is one other mesh, starting at the injective $P_4^4$ on the left
  which has no indecomposable end term.
\end{example}

\subsection{The adjoint isomorphism between covers and cokernels}

This adjoint isomorphism will be used to identify objects $(V,U)$ where the factor $V/U$
contains a summand of type $P^\ell$.

\medskip
In the submodule category $\mathcal S(R)$, let $\mathcal I$
be the categorical ideal of all maps which factor through an object of the form $(A=A)$
where $A\in\mod R$.  It is easy to see that
$$\mathcal I((U\subset V),(X\subset Y))\;=\;\{f\in\Hom_R(V,Y):f(V)\subset X\}.$$

\begin{lemma}
  A collection of assignments $M\mapsto(K\subset P)$ where
  $0\to K\to P\to M\to 0$ is a projective resolution for $M\in\mod R$
  defines a functor
  $$\Cov: \quad \mod R\to \mathcal S(R)/\mathcal I.$$
\end{lemma}

\begin{proof}
  Under $\Cov$, a map $h:M\to M'$ is sent to the class of the commuative square
  given by a lifting.  This assignment is well defined.

  \smallskip Namely, given a projective resolution $0\to K\to P\to M\to 0$,
  a short exact sequence $0\to U\to V\to W\to 0$, a map $h:M\to W$
  and two liftings $g_1,g_2:P\to V$, $f_1,f_2:K\to U$, then their difference
  $g=g_2-g_1$, $f=f_2-f_1$ gives rise to a commutative diagram
  $$
  \xymatrix{ 0 \ar[r] & K \ar[r] \ar[d]_f & P \ar[r] \ar[d]^g & M \ar[r] \ar[d]^0 & 0 \\
    0 \ar[r] & U \ar[r]^u & V \ar[r] & W \ar[r] & 0 }
  $$
  where the map $g$ factors over $u$.  Hence $(f,g)$, as a map in $\mathcal S(R)$,
  factors through the object $(U=U)$.
\end{proof}

\begin{lemma}
  \mylabel{lemma-I-exact}
  Let $0\to A\to B\to C\to 0$ be a short exact sequence in $\mod R$ and
  $$
  \xymatrix{ 0 \ar[r] & K_A \ar[r] \ar[d] & K_B \ar[r] \ar[d] & K_C \ar[r] \ar[d] & 0 \\
    0 \ar[r] & P_A \ar[r] &  P_A\oplus P_C \ar[r] & P_C \ar[r] & 0 }
  $$
  be the corresponding sequence in $\mathcal S(R)$, say obtained by using the
  Horseshoe Lemma. Then for any embedding $X:(U\subset V)$, the following sequence is split exact
  $$0\to\mathcal I((K_A\subset P_A),X)\to \mathcal I((K_B\subset P_A\oplus P_C),X)
  \to \mathcal I((K_C\subset P_C),X)\to 0.$$
\end{lemma}

\begin{proof}
  The last sequence simplifies to
  $$0\to \Hom_R(P_A,U)\to \Hom_R(P_A\oplus P_C,U)\to \Hom_R(P_C,U)\to 0.$$
\end{proof}

\begin{proposition}
  \mylabel{prop-cok-cov}
  The functors $R=\Cok, L=\Cov$ form an adjoint pair
  $$
  \xymatrix{\mathcal S(R)/\mathcal I\quad \ar@/^/[r]^{{\rm Cok}} & \quad \mod R \ar@/^/[l]^{{\rm Cov}} }
  $$
\end{proposition}

\begin{proof}
  Let $0\to K\to P\to M\to 0$ be a projective cover and $(U\subset V)$ an embedding
  with cokernel $W$.  We show that the map
  $$\alpha:\;\Hom_{\mathcal S/\mathcal I}(\Cov(M),(U\subset V))\;\to\;\Hom_R(M, \Cok(U\subset V))$$
  is an isomorphism.  The map $\alpha$ is defined since
  if a map $(f,g)$ in $\mathcal S(R)$ factors through an object
  of the form $(U,U)$, then the cokernel map is zero (since $\Cok:\mathcal S(R)\to \mod R$
  is a functor).

  \smallskip To see that $\alpha$ is invertible, 
  consider $h:M\to W$.   We have seen in the proof of the first lemma above
  that the difference of any two liftings factors through the embedding $(U,U)$.
  Hence if $(f,g)$ is a lifting for $h$ in the sense that the diagram is commutative
  $$
  \xymatrix{ 0 \ar[r] & K \ar[r] \ar[d]_f & P \ar[r] \ar[d]^g & M \ar[r] \ar[d]^h & 0 \\
    0 \ar[r] & U \ar[r] & V \ar[r] & W \ar[r] & 0 }
  $$
  then the morphism $(f,g)$ in $\mathcal S(R)/\mathcal I$ is determined uniquely.
  \end{proof}

\subsection{The multiplicity of $P^m$ as a summand of $V/U$}
\mylabel{sec-factor}

We return to the case where $R=k[T]/(T^n)$
and consider the Auslander-Reiten sequence (*) in $\mod R$ from Section~\ref{sec-ambient},
$$ 0\longrightarrow A\stackrel f\longrightarrow B\longrightarrow C\longrightarrow 0,$$
and an embedding $(U\subset V)$ with cokernel $W$.
We have seen that the multiplicity of $A$ as a direct summand of the $R$-module $W$
is $$\mu_A W=\dim\Cok\Hom_R(f,W).$$

\medskip
Using the Horseshoe Lemma as in the previous subsection, we obtain the sequence
$$\mathcal E: \quad (K_A\subset P_A)\stackrel{\tilde f}\longrightarrow
(K_B\subset P_A\oplus P_C)\longrightarrow (K_C\subset P_C)\longrightarrow 0;$$
we denote the objects in this sequence $\tilde A$, $\tilde B$, and $\tilde C$
and write $X=(U\subset V)$.

\medskip
From the adjoint isomorphisms in Proposition~\ref{prop-RA} and Proposition~\ref{prop-cok-cov}
and from Lemma~\ref{lemma-I-exact} we obtain
\begin{eqnarray*}\mu_A(W) & = &
  \dim\Cok\Hom_R(f,W)\\
  & = & \hom_R(C,W)-\hom_R(B,W)+\hom_R(A,W)\\
  & = & \hom_{\mathcal S/\mathcal I}(\tilde C, X)-\hom_{\mathcal S/\mathcal I}(\tilde B,X)
  +\hom_{\mathcal S/\mathcal I}(\tilde A,X)\\
  & = & \hom_{\mathcal S}(\tilde C,X)-\hom_{\mathcal S}(\tilde B,X)+\hom_{\mathcal S}(\tilde A,X)\\
  & = & \dim\Cok\Hom_{\mathcal S}(\tilde f, X).
\end{eqnarray*}
where we write $\hom$ for $\dim\Hom$.  This shows that the multiplicity of $A$ as a summand
of the factor is the covariant defect of the sequence $\mathcal E$ at the embedding $(U\subset V)$:
$$\mu_A(V/U) \;=\; e_*(U\subset V).$$

\begin{example}
  We return to the example in Section~\ref{sec-ambient}.
  In the case where $n=4$ and $m=3$,
  The Auslander-Reiten sequence in $\mathcal N(4)$ starting at $P^3$,
  $$0\to P^3\to P^2\oplus P^4\to P^3\to 0,$$
  gives rise to the short right exact sequence in $\mathcal S(4)$,
  $$\mathcal E:\quad P_1^4\stackrel{\tilde f}\longrightarrow P_2^4\oplus P_0^4\longrightarrow P_1^4\longrightarrow 0,$$
  hence for an object $X=(U\subset V)$, the multiplicity of $P^3$ as a summand of $V/U$ is $e_*(X)$.

  \smallskip
  The five objects $(V,U)$ marked with a large dashed semicircle have $P^3$ as a direct summand
  in the factor space $V/U$.  The above sequence
  $\mathcal E: P_1^4\to P_2^4\oplus P_0^4\to P_1^4\to 0$ gives rise to the hammock with source
  $X^f=P_1^4$, tangents at the meshes ending in $Y_1^f=P_2^4$ and $Y_2^f=P_0^4$ and sink at $\tau Z^f$ where
  $Z^f=P_1^4$.  There is an addtional tangent with no indecomposable end term at the injective object $P_0^4$
  labeled $Y_2^a$.

\bigskip
We would like to point out that the marked objects visualize a periodicity in the stable
part of the Auslander-Reiten quiver:  If the object $X$ is marked with a large solid semicircle,
then $\tau X$ has a large dotted semicircle, $\tau^2 X$ a small dashed semicircle,
$\tau^3 X$ a small solid semicircle, $\tau^4 X$ a small dotted semicircle, and $\tau^5 X$
a large dashed semicircle.  Finally, $\tau^6X\cong X$ is again marked with a large solid semicircle.
This illustrates \cite[Theorem~5.1]{rs2}:

\smallskip
Consider the indecomposable non-projective object $X=(U\subset V)$ as a triangle in the stable category $\overline{\mathcal N(n)}$
of the Frobenius category $\mathcal N(n)$,
$$X:\qquad U \to V\to V/U\to \Omega U\to $$
then the Auslander-Reiten translate for $X$, considered as a triangle in the stable category,
is given by rotating the above triangle:
$$\tau_{\mathcal S} X:\qquad V\to V/U\to \Omega U\to \Omega V\to $$
Thus, a non-projective summand $A$ of the ambient space in $X$ becomes a summand of the subspace of
$\tau_{\mathcal S} X$, then $\Omega_{\mathcal N}A$ becomes a summand of the factor space of $\tau^2_{\mathcal S} X$ etc.
\end{example}

\subsection{Example: A hammock with no isolated source or sink}
\mylabel{sec-no-isolated}

Consider the hammock function indicated by small numbers in the Auslander-Reiten quiver
studied before, it is pictured in Figure~\ref{figure-no-sink}.

\def\scale#1{\makebox[0pt][c]{\beginpicture\setcoordinatesystem units <1.8mm,1.8mm>#1
    \endpicture}}
\def\arqfoure{\beginpicture\setcoordinatesystem units <.45cm,.45cm>
  \multiput{} at -2 0  25 15 /
  \arqfourblack
  \multiput{$\green{1}$} at 4.7 3.7  16.7 8.5  4.6 6.5  16.6 6.5
          5 8.5  17 3.5  8.5 3.3  20.5 8.3  8.6 6.3  20.6 6.3  8.7 8.5  20.7 3.6  
          10.7 1.6  22.7 10.6  10.5 5.3  22.5 5.3  10.7 10.6  22.7 1.6
          .6 12.5  12.6 -.5  24.6 12.5  .6 8.5  12.6 3.5  24.6 8.5  .9 3.5  12.9 8.5  24.9 3.5
          2.7 10.7  14.7 1.7  2.5 5.3  14.5 5.3  2.7 1.7  14.7 10.7 /
  \multiput{$\green{2}$} at 6.9 5.6  18.9 5.6 /
  \multiput{\huge\bf\green{*}} at 0.1 2.4  12.1 10.4  24.1 2.4 /
  \endpicture
}

\begin{figure}[ht]
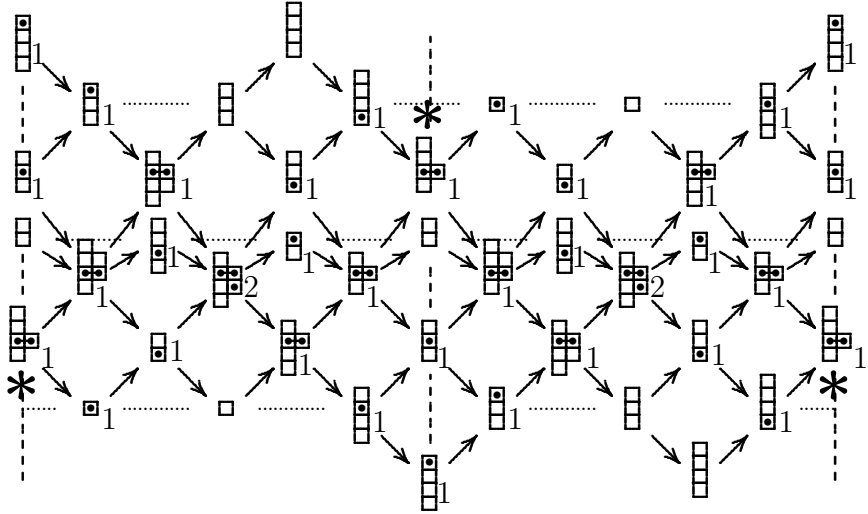

\begin{center}
  \arqfoure
\end{center}
\caption{Hammock function with no isolated source or sink}
\label{figure-no-sink}
\end{figure}

The function is given by counting for an object $(V,U)$ the number of indecomposable direct
summands of the subspace $U$.  This is the number of generators of $U$ as a $k[T]$-module
which are pictured as dots or connected sequences of dots, see the introduction in Section~\ref{sec-invariant}.

\medskip
Note that this hammock function is additive on each complete mesh except the one marked (*).
However, the hammock does not have an isolated source or an isolated sink.

\medskip
As the hammock function counts the number of indecomosable summands of $U$, it is given
by the short right exact sequence $\mathcal E: P_1^1\to 0\to 0\to 0$.
Hence the hammock has a non-isolated source at $P_1^1$, no tangent with an indecomposable end term,
and no non-injective sink.

\medskip
The picture in Figure~\ref{figure-no-sink2} shows this hammock as a dotted region, the source is labeled by $Z$.
Since it stretches over several copies
of the fundamental domain, there is an overlap.  We picture a shifted copy of the hammock,
encircled by a dotted line.  The union of the two regions yields the hammock function in the
above picture, note that the value is 2 where the two regions overlap.
The one object which occurs in both regions has been pictured as $P_{31}^{42}$ in the introduction to this section.
It is the unique indecomposable object
in $\mathcal S(4)$ for which the subspace as a $k[T]$-module has two indecomposable summands.

\def\scale#1{\makebox[0pt][c]{\beginpicture\setcoordinatesystem units <1.8mm,1.8mm>#1
    \endpicture}}
\def\arqfourf{\beginpicture\setcoordinatesystem units <.45cm,.45cm>
  \multiput{} at -2 0  25 15 /
  \arqfourblack
  \put{$Z$} at 3 1
  \put{$\tilde Z$} at 23.5 1
  \put{$\tilde Y$} at 19 -.5
  \setshadesymbol <z,z,z,z> ({$\red{.}$})
  \setshadegrid span <.7mm> point at 0 0
  \setlinear
  \vshade 1 1.5 2.5  1.5 1 3   2.5 1 4  5.5 4 7  7 5.5 9.5  7.5 6 10
  8.5 6  11  9.5 5 12  10.5 4 12  11 3.5 11.5  /
  \vshade 11 3.5 6.5  11.5 3 6  12.5 2 6  13 1.5 6.5 /
  \vshade 11 7.5 11.5  11.5 8 11  12.5 8 10  13 7.5 9.5 /
  \vshade 13 1.5 9.5  13.5 1 9  14.5 1 8  17.5 4 8  19 5.5 9.5  19.05 6.5 9.5
  19.5 8 10  20.5 8 11  23.5 11 14  24.5 12 14  25 12.5 13.5 /
  \setdots<2pt>
  \plot 12.8 11  13 11.5  13.5 12  14.5 12  15.5 11  16.5 10  17.5 9  18.5 8  20.5 8  21.5 8  22.5 9  23 9.5
  23.5 10  24.5 11  25 11.5 /
  \plot 12.8 11  13 10.5  13.5 10  14.5 9  15.5 8  16.5 8  17.5 5  18.5 4  20.5 2  21.5 1  22.5 1  23 1.5
  23.5 2  24.5 3  25 3.5 /
  \plot 25 7.5  24.5 8  23.5 8  23 7.5  22.8 7  23 6.5  23.5 6  24 6  25 6.5 /
  \plot -1 7.5  -.5 8  .5 8  1 7.5  1.2 7  1 6.5  .5 6  -.5 6  -1 6.5 /
  \plot -1 9.5  -.5 10  .5 11  1 11.5  1.5 12  2.5 12  3 11.3  5 9.5  11.5 2  12.5 1  13 .5  13.2 0
  13 -.5  12.5 -1  11.5 -1  5 5.5  3 5.5  2.5 5  1.5 4  1 3.5  .5 3  -.5 2  -1 1.5 /
  %
%
  \endpicture
}

\begin{figure}[ht]
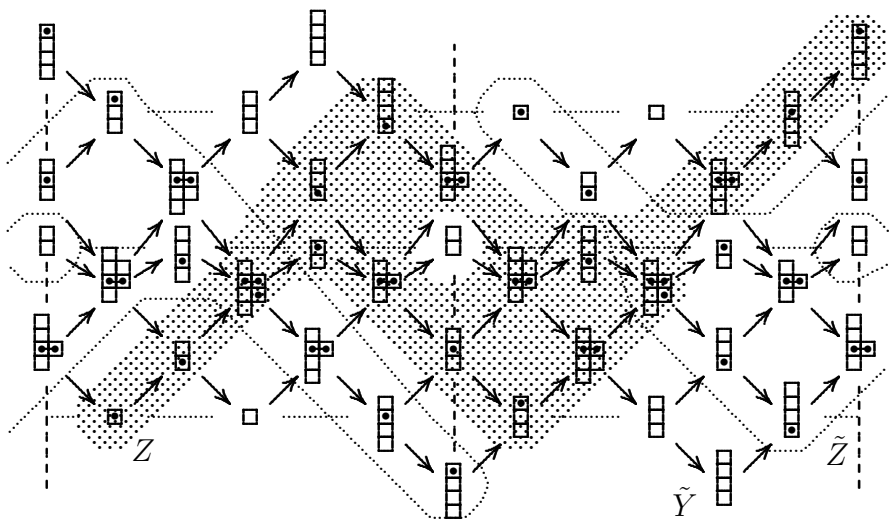

\begin{center}
  \arqfourf
\end{center}
\caption{Hammock counting the indecomposable summands of the subspace $U$}
\label{figure-no-sink2}
\end{figure}

\medskip
Note that the hammock is symmetric with respect to the vertical line through (*).
The same hammock function is obtained from the contravariant functor
given by the short left exact sequence
$\mathcal E: 0\to 0\to P_0^4\stackrel v\to P_1^4$.
It has a sink at $\tilde Z=P_1^4$, a tangent starting at the injective module
$\tilde Y=P_0^4$, and no non-projective source.


\bigskip
{\it Address of the author:}\\[1ex]
Department of Mathematical Sciences\\
Florida Atlantic University\\
777 Glades Road\\
Boca Raton, Florida 33431\\[1ex]
E-mail: {\tt markusschmidmeier@gmail.com}

\end{document}